\documentclass{amsart}
\usepackage{amssymb, amsmath, amsthm}

\usepackage[final, hypertex]{hyperref}

\newtheorem{theorem}{Theorem}[section]
\newtheorem{prop}[theorem]{Proposition}
\newtheorem{lemma}[theorem]{Lemma}
\newtheorem{cor}[theorem]{Corollary}

\theoremstyle{definition}
\newtheorem{dhef}[theorem]{Definition}
\newtheorem{example}[theorem]{Example}

\theoremstyle{remark}
\newtheorem{remark}[theorem]{Remark}

\numberwithin{equation}{section}

\def\R{\mathbb{R}}

\def\N{\mathbb{N}}
\def\U{\mathcal{U}}

\def\pscal#1#2{\left\langle#1,\,#2\right\rangle}

\def\haus{\mathcal{H}^{n-1}}
\def\dimh{\textrm{dim}_{\mathcal{H}}}

\def\convopen{\mathcal{K}^n_0}

\def\gau#1{\rho_{#1}}
\def\gauge{\rho}
\def\pgauge{\gauge^0}
\def\dists{\delta_S}
\def\distb#1{d_{#1}}
\def\dist{\distb{\Omega}}
\def\sdist{\dist^s}
\def\cut{m}
\def\uv{\nu}
\def\curv{\kappa}
\def\curvg{\tilde{\kappa}}
\def\nor{\nu}
\def\leno{\tilde{l}}
\def\len{l}

\def\ridge{\mathcal{R}}

\def\Lip{\textrm{Lip}}
\def\bw{\overline{W}}

\def\vf{v_f}
\def\meas{\mathcal{L}^n}
\def\matr{\mathcal{M}}

\DeclareMathOperator{\dive}{div}

\DeclareMathOperator{\diag}{diag}
\DeclareMathOperator{\proj}{\Pi}

\DeclareMathOperator{\cutl}{Cut}

\begin{document}
\title[Distance function]%
{The distance function from the boundary\\ in a Minkowski space}

\author[G.~Crasta]{Graziano Crasta}
\address{Dipartimento di Matematica ``G.\ Castelnuovo'', Univ.\ di Roma I\\
P.le A.\ Moro 2 -- 00185 Roma (Italy)}
\email[Graziano Crasta]{crasta@mat.uniroma1.it}

\author[A.~Malusa]{Annalisa Malusa}
\email[Annalisa Malusa]{malusa@mat.uniroma1.it}

\date{May 31, 2005}

\keywords{Distance function, Minkowski structure, cut locus,
Hamilton-Jacobi equations}
\subjclass[2000]{Primary 35A30;
Secondary 26B05, 32F45, 35C05, 49L25, 58J60}

\begin{abstract}
Let the space $\mathbb{R}^n$ be endowed with a Minkowski
structure $M$ (that is $M\colon \mathbb{R}^n \to [0,+\infty)$ is the gauge
function of a compact convex set having the origin as an interior
point, and with boundary of class $C^2$), and let $d^M(x,y)$ be
the (asymmetric) distance associated to $M$. Given an open domain
$\Omega\subset\mathbb{R}^n$ of class $C^2$, let $d_{\Omega}(x) :=
\inf\{d^M(x,y);\ y\in\partial\Omega\}$ be the Minkowski distance
of a point $x\in\Omega$ from the boundary of $\Omega$.
We prove that a suitable extension of $d_{\Omega}$ to $\mathbb{R}^n$ (which plays
the r\"ole of a signed Minkowski distance to $\partial \Omega$)
is of class $C^2$ in a tubular neighborhood of $\partial \Omega$,
and that $d_{\Omega}$ is of class $C^2$ outside the cut locus
of $\partial\Omega$
(that is the closure of the set of points of non--differentiability
of $d_{\Omega}$ in $\Omega$).
In addition, we prove that the cut locus of $\partial \Omega$
has Lebesgue measure zero, and that
$\Omega$ can be decomposed, up to this set of vanishing measure,
into geodesics starting from
$\partial\Omega$ and going into $\Omega$ along the
normal direction (with respect to the Minkowski distance).
We compute explicitly the Jacobian determinant
of the change of variables that associates to
every point $x\in \Omega$ outside the cut locus
the pair $(p(x), d_{\Omega}(x))$, where
$p(x)$ denotes the (unique) projection of $x$ on
$\partial\Omega$, and we apply
these techniques to the
analysis of PDEs of Monge-Kantorovich type
arising from problems in optimal transportation theory
and shape optimization.
\end{abstract}

\maketitle

\section{Introduction}

In recent years the study of the distance function from the
boundary has
attracted the attention of many researchers
coming from different areas of mathematical analysis.
We mention, among other papers, \cite{LN} in the framework
of general first order Hamilton--Jacobi equations and
Finsler geometry, \cite{IT} and \cite{MM} in the
case of Riemannian manifolds, \cite{CSW} and \cite{PRT}
for the point of view of non--smooth analysis
in Hilbert spaces,
\cite{EH} for results applied to the theory of Sobolev
spaces, and \cite{CFGH} for applications to causality
theory.

The results of the present paper are mainly motivated
by their applications to the analysis of PDEs.
Let $\Omega\subset\R^n$, $n\geq 2$, be a smooth open domain
(i.e., a nonempty open bounded connected subset of $\R^n$,
with sufficiently smooth boundary),
and let
$H\colon \overline{\Omega}\times\R^n\to\R$
be a smooth hamiltonian, such that
for every $x\in\overline{\Omega}$ the sublevel
$K(x) := \{p\in\R^n;\ H(x,p) \leq 1\}$
is a compact convex set
having the origin as an interior point,
with $C^2$ boundary and
with strictly positive principal curvatures at any point
($K(x)\in C^2_+$ for short).
It is well known that the function
\[
F(x, v) :=
\max\{\pscal{v}{p}; p\in K(x)\},
\qquad x\in\overline{\Omega}, \ v\in\R^n
\]
is a \textsl{Finsler structure} on $\overline{\Omega}$
(see \cite{BCS}), that is,
a smooth function defined in $\overline{\Omega}\times(\R^n\setminus\{0\})$,
positively $1$-homogeneous with respect to $v$,
and such that the Hessian matrix
$\partial^2_{v_i v_j} F^2$ is positive definite at every point
of $\overline{\Omega}\times(\R^n\setminus\{0\})$.

For every pair of points $x,y\in \overline{\Omega}$,
we can define the geodesic distance from $y$ to $x$
as
\[
L(x,y) :=
\inf\left\{
\int_0^1 F(\xi(t), \dot{\xi}(t))\, dt;\
\xi\in S_{x,y}\right\}\,,
\]
where $S_{x,y}$ is the set of all
Lipschitz curves $\xi\colon [0,1]\to \overline{\Omega}$
such that $\xi(0) = y$ and $\xi(1) = x$.
It is well known that the distance from the boundary of $\Omega$,
defined by
\begin{equation}\label{f:gendist}
\dist(x) := \inf_{y\in\partial\Omega} L(x,y),
\qquad x\in\overline{\Omega}
\end{equation}
is the unique viscosity solution of the Hamilton-Jacobi equation
\begin{equation}\label{f:HJi}
\begin{cases}
H(x, Du(x)) = 1
&\textrm{in $\Omega$},\\
u = 0
&\textrm{on $\partial\Omega$}
\end{cases}
\end{equation}
(see \cite{Li}).
Moreover,
for every $x\in\Omega$ the infimum in
(\ref{f:gendist}) is achieved at a point
$y\in\partial\Omega$ that can be joined
to $x$ by a geodesic going into $\Omega$
and starting from $y$ along the
``normal'' direction to $\partial\Omega$
(see \cite{LN}).

It is worth to remark that,
in the simple example $H(x,p) = |p|$,
then $F(x,v) = |v|$,
$(\R^n, F)$ is the standard Euclidean space,
and (\ref{f:HJi}) is the eikonal equation.
Hence $\dist$ is the Euclidean distance of $x$
from $\partial\Omega$, that is
\[
\dist(x) = \inf_{y\in\partial\Omega} |x-y|\,,
\qquad x\in\overline{\Omega}\,.
\]

Another case, that will be the one considered in our analysis,
concerns autonomous Hamiltonians of the form
$H(x,p) = \gauge(p)$,
where $\gauge(p) := \inf\{t\geq 0;\ p\in t K\}$
is the gauge function of a fixed
compact convex set $K\in C^2_+$.
In this case the function
$F(x,v)$ coincides with the gauge function $\pgauge(v)$
of the polar set $K^0$ of $K$, and
$(\R^n, F)$ is a Minkowski space
(see \cite[Chap.~14]{BCS}).
The function
\[
\dist(x) = \inf_{y\in\partial\Omega} \pgauge(x-y)\,,
\qquad x\in\overline{\Omega}\,,
\]
which is the unique viscosity solution of
\[
\begin{cases}
\gauge(Du)=1 & \textrm{in\ }\Omega\,, \\
u=0 & \textrm{on\ }\partial \Omega\,,
\end{cases}
\]
will be called the
Minkowski distance from $\partial\Omega$.
This non--symmetric distance function from the boundary
is exactly the object we shall deal with in this paper.

Our results can be divided into two main groups,
the first one devoted to the differentiability properties
of $\dist$, and the other one
on the regularity of the closure $\overline{\Sigma}$
of the set $\Sigma$ of those points in $\Omega$
where $\dist$ is not differentiable.

Concerning the differentiability of $\dist$,
we consider a signed distance $\sdist$, defined
in the whole $\R^n$, which extends $\dist$, and we prove that
$\sdist$ is of class $C^2$ in a tubular neighborhood
of $\partial\Omega$ (see Theorem~\ref{t:regd}).
This result provides a good definition of $D\dist$ and $D^2\dist$
on $\partial\Omega$.
For every $y\in\partial\Omega$, $D\dist(y)$ is proportional
to the inward normal $\nor(y)$ to $\partial \Omega$ at $y$
(see Lemma~\ref{l:pis}).
Moreover, due to the Minkowskian structure of the space, for every
$x\in\Omega$ and for every projection $x_0$ of $x$
on $\partial \Omega$ the
geodesic jointing $x_0$ with $x$ is a segment starting from
$x_0$ and going into $\Omega$ along the ``normal''
direction, which is $D\gauge(\nor(x_0))$ (see Proposition~\ref{p:pnc}).
Thanks to these regularity results, for every $y\in\partial\Omega$
we are able to define and study
a linear mapping acting on the tangent space to $\partial
\Omega$ at $y$ which plays the r\"ole of the Weingarten map
in Riemannian geometry, that is it has real eigenvalues which are
normal curvatures with respect to the Minkowski distance.

For what concerns the singular set,
a well known regularity result is
the rectifiability of $\Sigma$.
Namely, since $\dist$ is a locally semiconcave
function on $\Omega$ (see \cite{BaCD, CaSi}),
then $\Sigma$ is $C^2$-rectifiable (see \cite{Alb}),
that is
it can be covered by a countable family of
embedded $C^2$ manifolds of dimension $n-1$,
with the exception of a set of vanishing $\haus$ measure
(here $\mathcal{H}^s$ denotes the $s$-dimensional
Hausdorff measure).
In general, this is the best regularity result
that we can expect for $\Sigma$,
and remains valid even without any smoothness assumption
on $\partial\Omega$.

On the other hand, the set $\overline{\Sigma}$
may behave badly.
In \cite{MM} it is exhibited a domain $\Omega\subset\R^2$,
of class $C^{1,1}$, such that
$\overline{\Sigma}$ has positive Lebesgue measure.
(In the same paper the authors prove
a rectifiability result in a Riemannian setting.)
Li and Nirenberg in \cite{LN} have shown that,
if $\Omega$ is of class $C^{2,1}$
and $F$ is of class $C^{\infty}$,
then $\haus(\overline{\Sigma})$ is finite,
hence its Hausdorff dimension does not exceed $n-1$.
For related results in Finsler spaces see also
\cite[\S 6]{CaSi} and \cite{Me}.

Our result in this direction is the following
(see Theorem~\ref{t:regul}):
If $\Omega$ and $F(x,p) = \pgauge(p)$
are of class $C^{2,\alpha}$, for some $\alpha\in [0,1]$,
then the Hausdorff dimension
of $\overline{\Sigma}$ does not exceed $n-\alpha$.
In the case $\alpha = 0$ (that is, $\Omega$ and $\pgauge$
of class $C^2$),
we show that $\overline{\Sigma}$ has vanishing
Lebesgue measure,
and that $\dist$ is of class $C^2$ on
$\Omega\setminus\overline{\Sigma}$
(see Corollary~\ref{c:mn} and
Theorem~\ref{t:regd2}).
Finally, the set $\Omega\setminus \overline{\Sigma}$
can be decomposed
into geodesics (segments) starting from
$\partial\Omega$ and going into $\Omega$ along the
normal direction.
We compute explicitly the Jacobian determinant
of the change of variables which associates to
every point $x\in \Omega\setminus\overline{\Sigma}$
the pair $(p(x), \dist(x))$, where
$p(x)$ denotes the (unique) projection of $x$ on
$\partial\Omega$, and we are able to perform
a change of variables in
multiple integrals in $\Omega$ (see Theorem \ref{t:chvar}).

As an application of these results,
we consider the following system of PDEs of Monge-Kantorovich type:
\begin{equation}\label{f:syst1in}
\begin{cases}
-\dive(v\, D\gauge(Du)) = f
&\textrm{in $\Omega$},\\
\gauge(Du)\leq 1
&\textrm{in $\Omega$},\\
\gauge(Du) = 1
&\textrm{in $\{v>0\}$},
\end{cases}
\end{equation}
where the source $f\geq 0$ is a continuous function in $\Omega$,
complemented with the conditions
\begin{equation}\label{f:syst2in}
\begin{cases}
u\geq 0,\
v\geq 0
&\textrm{in $\Omega$},\\
u=0
&\textrm{on $\partial\Omega$}.
\end{cases}
\end{equation}
The first equation in (\ref{f:syst1in}) has to be
understood in the sense of distributions, whereas
$u$ is a viscosity solution to the Hamilton--Jacobi
equation
$\gauge(Du) = 1$ in the set $\{v>0\}$.
We look for a solution $(u,v)$ to
(\ref{f:syst1in})--(\ref{f:syst2in}) in the class of continuous
and non-negative functions.

This system of PDEs
arises in problems
of shape optimization (see \cite{BoBu}).
In the case $\gauge(\xi) = |\xi|$,
(\ref{f:syst1in})--(\ref{f:syst2in}) describes the stationary
solutions of models in granular matter
theory (see \cite{CCCG}).

In section~\ref{s:PDE}, using a change of variable formula, we shall
explicitly construct a non-negative and continuous function
$\vf$ (defined in (\ref{f:vf})) such that the pair $(\dist, \vf)$
is a solution to (\ref{f:syst1in})--(\ref{f:syst2in}),
extending a result proved in \cite{CCCG}
in the case $\gauge(\xi) = |\xi|$.

\section{Notation and Preliminaries}

\subsection{Basic notation}
The standard scalar product of two
vectors $x,y\in\R^n$
is denoted by $\pscal{x}{y}$,
and $|x|$ denotes the
Euclidean norm of $x\in\R^n$.
By $S^{n-1}$ we denote the set of unit vectors
of $\R^n$, and by $\matr_k$ the set of
$k\times k$ square matrices.
We shall denote by $(e_1,\ldots,e_n)$
the standard basis of $\R^n$.
Given two points $x,y\in\R^n$,
$[x,y]$ will denote the closed segment joining $x$ to $y$,
while $(x,y)$ will denote the same segment without
the endpoints.
If $A$, $B\subset\R^n$, $x\in\R^n$ and $t\in\R$, we define
$A+x = \{a+x;\ a\in A\}$,
$t\, A = \{t\, a;\ a\in A\}$ and
$A+B = \{a+b;\ a\in A,\ b\in B\}$.

As is customary, $B_r(x_0)$ and $\overline{B}_r(x_0)$
are respectively the open and the closed ball
centered at $x_0$ and with radius $r>0$.
Given two vectors $v,w\in\R^n$, the symbol
$v\otimes w$ will denote their tensor product,
that is, the linear application from
$\R^n$ to $\R^n$ defined by
$(v\otimes w)(x) = v\,\pscal{w}{x}$.

We shall denote by $\meas(A)$ and $\mathcal{H}^s(A)$,
respectively, the Lebesgue measure and the
$s$-dimensional Hausdorff measure of a set
$A\subset\R^n$.

Given $A\subset\R^n$,
we shall denote by $\Lip(A)$, $C(A)$ and $C^k(A)$, $k\in\N$,
the sets of functions $u\colon A\to\R$
that are respectively Lipschitz continuous, continuous
and $k$-times continuously differentiable in $A$.
Moreover, $C^{\infty}(A)$ will denote the set of functions
of class $C^k(A)$ for every $k\in\N$,
while $C^{k,\alpha}(A)$ will be the set of functions
of class $C^k(A)$ with H\"older continuous $k$-th partial derivatives
with exponent $\alpha\in [0,1]$.

A bounded open set $A\subset\R^n$
(or, equivalently, its closure
$\overline{A}$ or its boundary $\partial A$)
is of class $C^k$, $k\in\N$,
if for every point $x_0\in\partial A$
there exists a ball $B=B_r(x_0)$ and a one-to-one
mapping $\psi\colon B\to D$ such that
$\psi\in C^k(B)$, $\psi^{-1}\in C^k(D)$,
$\psi(B\cap A)\subseteq\{x\in\R^n;\ x_n > 0\}$,
$\psi(B\cap\partial A)\subseteq\{x\in\R^n;\ x_n = 0\}$.
If the maps $\psi$ and $\psi^{-1}$ are of class
$C^{\infty}$ or $C^{k,\alpha}$ ($k\in\N$, $\alpha\in [0,1]$),
then $A$ is said to be of class
$C^{\infty}$ or $C^{k,\alpha}$ respectively.

\subsection{Differential geometry}
We recall briefly some elementary facts
from differential geometry of hypersurfaces
of class $C^2$
(see e.g.~\cite{Th}).
Let $A\subset\R^n$ be a bounded open set of class $C^2$.
For every $x\in\partial A$, we denote respectively by $\nor(x)$
and $T_x A$
the unique inward unit normal vector
and the tangent space
of $\partial A$ at $x$.
The map $\nor\colon\partial A\to S^{n-1}$
is called the
\textsl{spherical image map}
(or \textsl{Gauss map}).
It is of class $C^1$ and,
for every $x\in\partial A$,
its differential $d\nor_x$
maps the tangent space $T_x A$ into itself.
The linear map
$L_x:= - d\nor_x\colon T_x A\to T_x A$ is called the
\textsl{Weingarten map}.
The bilinear form defined on $T_x A$
by
$S_x(v,w) = \pscal{L_x\, v}{w}$,
$v,w\in T_x A$,
is the
\textsl{second fundamental form} of $\partial A$ at $x$.
The geometric meaning of the Weingarten map
is the following.
For every $v\in T_x A$ with unit norm,
$S_x(v,v)$ is equal to the normal curvature of
$\partial A$ at $x$ in the direction $v$,
that is,
$S_x(v,v) = \pscal{\ddot{\xi}(0)}{\nor(x)}$,
where $\xi(t)$ is any parameterized curve in
$\partial A$ such that $\xi(0) = x$
and $\dot{\xi}(0) = v$.
The eigenvalues $\curv_1(x),\ldots,\curv_{n-1}(x)$
of the Weingarten map $L_x$ are,
by definition, the
\textsl{principal curvatures} of $\partial A$ at $x$.
The corresponding eigenvectors are called
the \textsl{principal directions}
of $\partial A$ at $x$.
It is readily shown that every $\curv_i(x)$ is
the normal curvature of
$\partial A$ at $x$ in the direction
of the corresponding eigenvector.
{}From the $C^2$ regularity assumption on the manifold $\partial A$,
it follows that the principal curvatures of $\partial A$
are continuous functions on $\partial A$.

\subsection{Convex geometry}
By $\convopen$ we denote the class of
nonempty, compact, convex subsets of $\R^n$
with the origin as an interior point.
We shall briefly refer to the elements of $\convopen$
as \textsl{convex bodies}.
The polar body of a convex body $K\in\convopen$
is defined by
\[
K^0 = \{p\in\R^n;\ \pscal{p}{x}\leq 1\ \forall x\in K\}\,.
\]
We recall that, if $K\in\convopen$, then
$K^0\in\convopen$ and $K^{00} = (K^0)^0 = K$
(see \cite[Thm.~1.6.1]{Sch}).

Given $K\in\convopen$ we define
its gauge function as
\[
\gau{K}(\xi) = \inf\{ t\geq 0;\ \xi\in t K\}\,.
\]
It is easily seen that
\[
\gau{K^0}(\xi) =
\sup\left\{\pscal{\xi}{p};\ p\in K\right\}\,,
\]
that is, the gauge function of the polar set $K^0$
coincides with the support function of the set $K$.
Let $0<c_1\leq c_2$ be such that
$\overline{B}_{c_2^{-1}}(0)\subseteq K \subseteq \overline{B}_{c_1^{-1}}(0)$.
Upon observing that $\xi/\gau{K}(\xi)\in K$ for every $\xi\neq 0$,
we get
\begin{equation}\label{f:brho}
c_1 {|\xi|}\leq\gau{K}(\xi)\leq c_2 {|\xi|}\,,\quad
\forall\xi\in\R^n.
\end{equation}

We say that $K\in\convopen$ is of class $C^2_+$
if $\partial K$ is of class $C^2$
and all the principal curvatures are strictly positive
functions on $\partial K$.
In this case, we define the $i$-th principal radius of
curvature at $x\in\partial K$ as
the reciprocal of the $i$-th principal curvature of $\partial K$
at $x$.
We remark that,
if $K$ is of class $C^2_+$, then $K^0$ is also of
class $C^2_+$
(see \cite[p.~111]{Sch}). Moreover, a convex body of class
$C^2_+$ is necessarily a strictly convex set.

Throughout the paper we assume that
\begin{equation}\label{f:ipoK}
K\in\convopen\
\textrm{is of class}\ C^2_+\,.
\end{equation}
Since $K$ will be kept fixed, from now on we shall use
the notation $\gauge = \gauge_K$ and $\pgauge=\gauge_{K^0}$.

We collect here some known properties
of $\gauge$ and $\pgauge$
that will be frequently used
in the sequel.

\begin{theorem}\label{t:sch}
Let $K$ satisfy $(\ref{f:ipoK})$. Then the following hold:
\par\noindent (i)
The functions $\gauge$ and $\pgauge$ are
convex, positively $1$-homoge\-neous
in $\R^n$, and
of class $C^2$ in $\R^n\setminus\{0\}$.
As a consequence,
\begin{gather*}
\gauge(t\, \xi) = t\, \gauge(\xi),\quad
D\gauge(t\, \xi) = D \gauge(\xi),\quad
D^2\gauge(t\, \xi) = \frac{1}{t}\, D^2\gauge(\xi),\\
\pgauge(t\, \xi) = t\, \pgauge(\xi),\
D\pgauge(t\, \xi) = D \pgauge(\xi),\
D^2\pgauge(t\, \xi) = \frac{1}{t}\, D^2\pgauge(\xi),
\end{gather*}
for every $\xi\in\R^n\setminus\{0\}$ and $t>0$.
Moreover
\begin{gather}
\pscal{D\gauge(\xi)}{\xi} = \gauge(\xi),\quad
\pscal{D\pgauge(\xi)}{\xi} = \pgauge(\xi)\,,
\quad\forall\xi\in\R^n\setminus\{0\}
\label{f:euler}\\
\label{f:eigenzero}
D^2\gauge(\xi)\, \xi = 0\,,\quad
D^2\pgauge(\xi)\, \xi = 0\,,
\qquad\forall \xi\in\R^n\setminus\{0\}\,.
\end{gather}
\par\noindent (ii)
For every $\xi,\eta\in\R^n$, we have
\begin{equation}\label{f:sublin}
\gauge(\xi+\eta)\leq \gauge(\xi) + \gauge(\eta),\quad
\pgauge(\xi+\eta)\leq \pgauge(\xi) + \pgauge(\eta),
\end{equation}
and equality holds if and only if $\xi$ and $\eta$ belong
to the same ray, that is,
$\xi = \lambda\, \eta$ or $\eta = \lambda\, \xi$
for some $\lambda\geq 0$.
\par\noindent (iii)
For every $\xi\neq 0$,
$D\gauge(\xi)$ belongs to $\partial K^0$, while $D\pgauge(\xi)$
belongs to $\partial K$. More precisely, $D\gauge(\xi)$
is the unique point of $\partial K^0$ such that
\[
\pscal{D\gauge(\xi)}{\xi} = \gauge(\xi),\ \textrm{and\ }
\pscal{x}{\xi} < \gauge(\xi)\ \forall x\in K^0,\ x\neq D\gauge(\xi)\,.
\]
Symmetrically,
the gradient of $D\pgauge(\xi)$ is the
unique point of $\partial K$ such that
\[
\pscal{D\pgauge(\xi)}{\xi} = \pgauge(\xi),\ \textrm{and\ }
\pscal{x}{\xi} < \pgauge(\xi)\ \forall x\in K,\ x\neq D\pgauge(\xi)\,.
\]
\par\noindent (iv)
The ei\-gen\-val\-ues of the second differential
$D^2\gauge$ at $\uv\in S^{n-1}$
are $0$, with corresponding eigenvector $\uv$,
and the principal radii of curvature of
$\partial K^0$ at the unique point $p\in\partial K^0$
at which $\uv$ is attained as an outward normal vector.
Symmetrically, the ei\-gen\-val\-ues of
$D^2\pgauge$ at $\uv\in S^{n-1}$
are $0$, with corresponding eigenvector $\uv$,
and the principal radii of curvature of
$\partial K$ at the unique point $p\in\partial K$
at which $\uv$ is attained as an outward normal vector.
\end{theorem}

\begin{proof}
(i) The convexity and the positive $1$-homogeneity are a consequence
of the definition.
The $C^2$ regularity in $\R^n\setminus\{0\}$
is proved in \cite[p.~106]{Sch}.
The identities follow upon observing that
$D\gauge$ and $D\pgauge$ are positively
$0$-homogeneous,
whereas $D^2\gauge$ and $D^2\pgauge$ are positively
$(-1)$-homogeneous.
The identity (\ref{f:euler}) is Euler's formula for positively
$1$-homogeneous functions,
whereas (\ref{f:eigenzero})
can be obtained differentiating (\ref{f:euler}).

\noindent
(ii)
Follows from convexity of $\gauge$ and $\pgauge$
and strict convexity of $K$ and $K^0$.

\noindent
(iii-iv)
See \cite{Sch},
Corollaries~1.7.3 and~2.5.2.
\end{proof}

\begin{lemma}\label{l:polare}
Let $K$ satisfy $(\ref{f:ipoK})$.
Then the identities
\begin{equation}\label{f:polare}
D\pgauge(D\gauge(\xi)) = \frac{\xi}{\gauge(\xi)}\,,\qquad
D\gauge(D\pgauge(\xi)) = \frac{\xi}{\pgauge(\xi)}
\end{equation}
hold for every $\xi\in\R^n\setminus\{0\}$.
\end{lemma}

\begin{proof}
Fixed $\xi\in\R^n\setminus \{0\}$, by Theorem~\ref{t:sch}(iii)
we have that $D\gauge(\xi)\in \partial K^0$, $D\pgauge(D \gauge(\xi))
\in \partial K$, hence $\pgauge(D\gauge(\xi))=1$ and
\begin{equation}\label{f:pol2}
\begin{split}
&\pscal{D\pgauge(D\gauge(\xi))}{D\gauge(\xi)} = 1\, \\
&\pscal{x}{D \gauge(\xi)} <1 \quad \forall x\in K,\ x\neq D\pgauge(D\gauge(\xi))\,.
\end{split}
\end{equation}
On the other hand by (\ref{f:euler}), the point $\xi/\gauge(\xi)\in \partial K$ satisfies
\[
\pscal{D\gauge(\xi)}{\frac{\xi}{\gauge(\xi)}}=1,
\]
which, together with (\ref{f:pol2}), implies the first identity in (\ref{f:polare}).
The second identity can be obtained from the first
interchanging the role of $K^0$ and $K$.
\end{proof}

\begin{remark}
It can be checked that Lemma~\ref{l:polare}
holds under the weaker assumption
that the sets $K,K^0\in\convopen$ be both strictly convex.
This is equivalent to require that
$K$ is a strictly convex body
of class $C^1$.
\end{remark}

\subsection{Nonsmooth analysis}
Let us recall the notions
of semiconcave function and
of viscosity solution
of Hamilton-Jacobi equations
(see for example \cite{BaCD, CaSi, Li}).
A function $u\colon A\to\R$ is said to be semiconcave
if there exists a constant $C>0$ such that
\[
\lambda u(x) + (1-\lambda) u(y)\leq
u(\lambda x + (1-\lambda)y) +
C\lambda(1-\lambda) |x-y|^2
\]
for every $\lambda\in [0,1]$ and every pair $x,y\in A$ such
that the segment $[x,y]$ is contained in $A$.
It is easy to check that this amounts to the concavity of the
map $x\mapsto u(x)-C|x|^2$ on every convex subset of $A$.
As a consequence, if $u$ is semiconcave in $A$
then it is locally Lipschitz continuous in $A$.

Let $A\subset\R^n$ be an open set, and
let $u\colon A\to\R$ be a continuous function.
The superdifferential $D^+u(x)$ and the subdifferential $D^-u(x)$
of $u$ at $x\in A$ are defined by
\[
\begin{split}
D^+u(x) & = \left\{
p\in\R^n;\ \limsup_{y\to x,\ y\in A}
\frac{u(y)-u(x)-\pscal{p}{y-x}}{|x-y|}\leq 0
\right\}\,,\\
D^-u(x) & = \left\{
p\in\R^n;\ \liminf_{y\to x,\ y\in A}
\frac{u(y)-u(x)-\pscal{p}{y-x}}{|x-y|}\geq 0
\right\}\,.
\end{split}
\]
If $u$ is differentiable at $x$,
then $D^+u(x) = D^- u(x) = \{Du(x)\}$
(see \cite[Lemma~1.8]{BaCD}).

Let $H\colon A\times\R\times\R^n\to\R$ be a
continuous function.
We say that $u\in C(A)$ is a \textsl{viscosity solution}
of
\begin{equation}\label{f:HJgen}
H(x, u, Du) = 0\qquad
\textrm{in}\ A,
\end{equation}
if, for every $x\in A$,
we have that $H(x,u(x),p)\leq 0$
for every $p\in D^+ u(x)$ and
$H(x,u(x),p)\geq 0$ for every $p\in D^- u(x)$.

\begin{example}
Given an open bounded set $A\subset\R^n$,
the Euclidean distance from the boundary of $A$, defined by
\[
u(x) = \inf_{y\in\partial A} |x-y|\,
\qquad x\in A,
\]
is the unique viscosity solution of the eikonal
equation $|Du| = 1$ in $A$ satisfying the
boundary condition $u=0$ on $\partial A$
(see \cite{BaCD, Li}).
\end{example}

Let $u\colon A\to\R$ be a locally Lipschitz function.
A vector $p\in\R^n$ is a \textsl{reachable gradient}
of $u$ at $x\in A$ if there exists a sequence
$(x_k)_k$ in $A\setminus\{x\}$ converging to $x$,
such that $u$ is differentiable at $x_k$
for every $k$ and $(Du(x_k))_k$ converges to $p$.
The set of all reachable gradients of $u$ at $x$
is denoted by $D^*u(x)$.

Since $u$ is locally Lipschitz, it is easily seen
that $D^*u(x)$ is compact for every $x\in A$.
Moreover, since $u$ is differentiable almost everywhere
in $A$, $D^* u(x)$ is nonempty for every $x\in A$.

\section{Minkowski distance from the boundary}

In the first part of this section we shall define
the Minkowski distance from a closed
set $S\subset\R^n$ and we shall prove some basic properties.
Then we shall specialize to the case $S=\partial\Omega$,
where $\Omega\subset\R^n$ is a nonempty bounded open set,
and we shall study some properties of the Minkowski
distance from $\partial\Omega$ in $\Omega$.

We recall that $K\in\convopen$ is a fixed convex body
of class $C^2_+$ with polar set $K^0$, and that $\gauge$,
$\pgauge$ are respectively the gauge function of $K$
and $K^0$.

Throughout this section, $S$ will be a
nonempty closed subset of $\R^n$.

\begin{dhef}[Minkowski distance from a set]
The Minkowski distance from $S$
is the function
$\dists\colon\R^n\to\R$ defined by
\begin{equation}\label{f:ds}
\dists(x) = \min_{y\in S} \pgauge(x-y)\,,
\qquad x\in\R^n\,.
\end{equation}
\end{dhef}

Some comments are in order.
Let us consider
the Euclidean distance from $S$,
defined by
\begin{equation}\label{f:euclid}
d^E_S(x) = \min_{y\in S} |x-y|\,,
\qquad x\in \R^n.
\end{equation}
It is clear that $d^E_S(x)$ measures the
(Euclidean) length of the shortest segment $[y,x]$
joining $x$ to a point $y\in S$.
In the definition of $\dists$,
the Euclidean distance $d^E(x,y) = |x-y|$ is replaced by
the Minkowski distance
$d^M(x,y) = \pgauge(x-y)$.
Besides the name,
the function $d^M$ is \textsl{not}
a distance in the usual sense.
Namely, it satisfies
\begin{gather*}
d^M(x,y)\geq 0, \ \textrm{and\ } d^M(x,y) = 0\ \textrm{only if \ }x=y \\
d^M(x,y)\leq d^M(x,z) + d^M(z,y)
\end{gather*}
for every $x,y,z\in\R^n$,
but in general it is not symmetric
(unless $K$ is a convex body symmetric
with respect to the origin).
Hence, in the Minkowski distance the ``length''
of a segment $[y,x]$, defined by
$d^M(x,y) = \gauge^0(x-y)$,
does not necessarily coincide with the length
of the segment $[x,y]$.
The space $(\R^n, d^M)$ is said to be a
Minkowski space.
It is well known that geodesics in Minkowskian spaces
are straight lines, so that these spaces are
geodesically both forward and backward complete
(see \cite[Chap.~14]{BCS}).

\begin{remark}
As we have underlined in the introduction,
we are interested in applications to nonlinear
PDEs, where the convex body $K$
plays the role of a constraint on the gradients
of the admissible functions.
More precisely, one often deals with the set of
Lipschitz continuous functions
$u\colon A\subset\R^n\to \R$ such that
$Du(x)\in K$ for almost
every $x\in A$, that is
$\gauge(Du(x))\leq 1$ for almost every $x\in A$.
We shall see that the
Minkowski distance from $S$ satisfies
$\gauge(D\dists(x)) = 1$ for
almost every $x\in A := \R^n\setminus S$
(see Proposition~\ref{p:basic}(i), and Theorem~\ref{t:sch}(iii)).
For this reason we prefer to define $\dists$
in terms of $\pgauge$ (and hence of $K^0$)
instead of using $\gauge$.
\end{remark}

Since $S$ is nonempty and closed,
for every $x\in\R^n$ the set
of projections of $x$ in $S$, defined by
\[
\proj_S(x) = \{y\in S;\ \dists(x) = \pgauge(x-y)\}\,,
\]
is nonempty and compact.

In the following proposition
we gather some basic properties of the
Min\-kow\-ski distance $\dists$.

\begin{prop}\label{p:basic}
The Minkowski distance $\dists$ from $S$
is a continuous function in $\R^n$,
locally semiconcave in $\R^n\setminus S$ and,
for every $x\not\in S$,
the following hold.
\begin{itemize}
\item[(i)]
$\dists$ is differentiable at $x$ if and only if
$\proj_S(x) = \{y\}$ is a singleton;
in this case
$D\dists(x) = D\pgauge(x-y)$ and
\[
D\gauge(D\dists(x)) = \frac{x-y}{\pgauge(x-y)}
= \frac{x-y}{\dists(x)}\,.
\]
\item[(ii)]
If $\proj_S(x)$ is not a singleton, then
$D^+ \dists(x)$ is the convex hull of the set
$\{D\pgauge(x-y);\ y\in\proj_S(x)\}$,
while $D^-\dists(x) = \emptyset$.
\item[(iii)]
For any $y\in\proj_S(x)$, and for any point
$z$ in the segment $(x,y)$ we have $\proj_S(z) = \{y\}$.
Hence $\dists$ is differentiable at
every point $z\in (x,y)$ and
$D\dists(z) = D\pgauge(x-y)$.
\item[(iv)]
$D^*\dists(x) = \{D\pgauge(x-y);\ y\in\proj_S(x)\}$.
\end{itemize}
\end{prop}

\begin{proof}
The proof below is a straightforward adaptation of
Proposition~2.14 in \cite{BaCD}
and Corollary~3.4.5 in \cite{CaSi},
where the case $K = \overline{B}_1(0)$ is considered.

The continuity of $\dists$ follows from the inequalities
\[
-\pgauge(z-x)\leq \dists(x) - \dists(z)\leq \pgauge(x-z),
\qquad
\forall x,z,\in\R^n\,.
\]
The other properties follow from the general theory concerning
the marginal functions
(see \cite[Proposition~2.13]{BaCD},
\cite[Proposition~3.4.4]{CaSi}).
Let $A=\R^n\setminus S$ and
$F(y,x) = \pgauge(x-y)$,
so that $\dists(x) = \min_{y\in S} F(y,x)$ for
every $x\in A$.
Since $\pgauge$ is of class $C^2$ in $\R^n\setminus\{0\}$,
it is clear that
$F(y,x)$, $D_x F(y,x) = D\pgauge(x-y)$ and
$D_{xx} F(y,x) = D^2\pgauge(x-y)$
are continuous functions in $S\times A$.
The local semiconcavity of $\dists$
in $\R^n\setminus S$ is a consequence of
the same property for marginal functions
proved in \cite[Proposition~3.4.1]{CaSi}.

If we define the sets
\begin{gather*}
M(x) := \{y\in S;\ \dists(x) = F(y,x)\} = \proj_S(x)\,,\\
Y(x) := \{D_x F(y,x);\ y\in M(x)\}
= \{D\pgauge(x-y);\ y\in\proj_S(x)\}\,,
\end{gather*}
properties (i) and (ii) are an immediate consequence
of \cite[Proposition~3.4.4]{CaSi} and Lemma~\ref{l:polare}.

In order to prove (iii),
let us start by proving that $y\in\proj_S(z)$ for every $z\in(x,y)$.
Namely, if not, there exists $y'\in S$ such that
$\pgauge(z-y)>\pgauge(z-y')$.
Since $z\in(x,y)$,
by Theorem~\ref{t:sch}(ii) we obtain
\[
\pgauge(x-y)=\pgauge(x-z)+\pgauge(z-y)>\pgauge(x-z)+\pgauge(z-y')
\geq \pgauge(x-y')\,,
\]
which contradicts the fact that $y\in\proj_S(x)$.
The same argument shows that $\proj_S(z)=\{y\}$.
Namely, if we suppose that there exists
$y'\in\proj_S(z)$, $y'\neq y$, we have that $x-z$
and $z-y'$ are not proportional, hence,
again from Theorem~\ref{t:sch}(ii), we get
\[
\pgauge(x-y') < \pgauge(x-z)+\pgauge(z-y')
= \pgauge(x-z) + \pgauge(z-y)
= \pgauge(x-y),
\]
which contradicts $y\in\proj_S(x)$.
Now the differentiability of $\dists$ at $z$ and the fact that
$D \dists(z)=D\pgauge(x-y)$ follow from (i) and from the
0--homogeneity of $D \pgauge$.

Finally, in order to prove (iv),
we have to show that $D^*\dists(x) = Y(x)$.
The inclusion $D^*\dists(x)\subset Y(x)$
follows from \cite[Proposition~3.4.4]{CaSi}.
Now let $p\in Y(x)$, and let us prove that $p\in D^*\dist(x)$.
By definition of $Y(x)$, there exists $y\in\proj_S(x)$
such that $p = D\pgauge(x-y)$.
{}From (iii) it we have that $\dists$ is differentiable
at any point $z\in (x,y)$, and $\dists(z) = p$,
hence $p$ is a reachable gradient at $x$.
\end{proof}

The following lemma, which states that the multifunction
$x \mapsto \proj_S(x)$ is sequentially upper semicontinuous,
is a straightforward
consequence of the continuity of $\pgauge$
and $\dists$.

\begin{lemma}\label{l:pusc}
If $(x_k)_k\subset\R^n$
is a sequence converging to a point $x\in\R^n$,
if $y_k\in\proj_S(x_k)$ for every $k$,
and if the sequence $(y_k)_k$ converges to $y$,
then $y\in\proj_S(x)$.
\end{lemma}

\begin{prop}\label{p:raylength}
Let $x\not\in S$.
For $y\in\proj_S(x)$ let us define
\begin{equation}\label{f:lenoy}
\leno_y(x) = {\dists(x)}\cdot
\sup\{t>0;\ y\in\proj_S(y+t(x-y))\}\,
\in (0,+\infty]\,,
\end{equation}
and, if $\leno_y(x) < +\infty$, the point
\begin{equation}\label{f:cuty}
\cut_y(x) = y + \leno_y(x)\, \frac{x-y}{\dists(x)}\, .
\end{equation}
Then $\leno_y(x)$ and $\cut_y(x)$ (when $\leno_y(x)$ is finite)
do not depend on
the choice of $y\in\proj_S(x)$.
\end{prop}

\begin{proof}
Observe that, if $z=y+t(x-y)$, $t>0$,
and $y\in\proj_S(z)$, then
\[
\dists(z) = \pgauge(z-y) = t\, \pgauge(x-y) = t\, \dists(x)\,.
\]
From the very definition
of $\leno_y(x)$
and Proposition~\ref{p:basic}(iii),
we have that $\leno_y(x)\geq \dists(x)$,
for every $y\in\proj_S(x)$.
We have two possibilities:
\begin{itemize}
\item[(i)]
$\leno_y(x) = \dists(x)$
for every $y\in\proj_S(x)$;
\item[(ii)]
$\leno_y(x) > \dists(x)$
for some $y\in\proj_S(x)$.
\end{itemize}
In case (i) it is plain that $\leno_y(x)$ does
not depend on the choice of $y\in\proj_S(x)$.
Furthermore, the right hand side of (\ref{f:cuty})
coincides with $x$ for every $y\in\proj_S(x)$,
hence $\cut_y(x)=x$ is also independent of the
choice of $y\in\proj_S(x)$.
Let us consider case (ii).
Let $y\in\proj_S(x)$ such that $\leno_y(x) > \dists(x)$.
By definition of $\leno_y(x)$,
there exists $t > 1$ such that,
setting $z=y+t(x-y)$,
$y\in\proj_S(z)$.
Since $x\in(y,z)$,
from Proposition~\ref{p:basic}(iii)
we have that $\proj_S(x) = \{y\}$
is a singleton.
This concludes the proof of the proposition.
\end{proof}

Thanks to Proposition~\ref{p:raylength},
the function
$\leno\colon\R^n\setminus S\to (0,+\infty]$,
\begin{equation}\label{f:leno}
\leno(x) := {\dists(x)}\cdot
\sup\{t>0;\ y\in\proj_S(y+t(x-y))\}\,
\qquad (y\in\proj_S(x))
\end{equation}
is well-defined.
The same holds true for
the function
$\cut\colon\{x\in\R^n\setminus S;\
\leno(x) < +\infty\}\to\Omega$
defined by
\begin{equation}\label{f:cut}
\cut(x) := y + \leno(x)\, \frac{x-y}{\dists(x)}\,
\qquad (y\in\proj_S(x))\,.
\end{equation}

\begin{remark}\label{r:leno}
It is straightforward to see that, given $x\not\in S$
and $y\in\proj_S(x)$,
\begin{equation}\label{f:leno2}
\leno(x) =
\sup\left\{\lambda>0;\
y\in\proj_S\left(x_{\lambda}
\right)\right\}\,,\quad
x_{\lambda} := y+\lambda\, \frac{x-y}{\pgauge(x-y)}\,.
\end{equation}
It is easily checked that $x_{\lambda} = x$
for $\lambda = \dists(x)$,
hence,
from Proposition~\ref{p:basic}(iii),
$\proj_S(x_{\lambda}) = \{y\}$
for every $\lambda\in [0,\dists(x))$.
As a consequence, we recover the inequality
$\leno(x)\geq\dists(x)$.
\end{remark}


\medskip
From now on we shall restrict our attention to
the case $S=\R^n\setminus\Omega$,
where $\Omega\subset\R^n$
will be a fixed nonempty connected bounded open set.

\begin{dhef}[Minkowski distance from the boundary]\label{d:dist}
\ The Min\-kow\-ski distance from the
boundary of $\Omega$ is defined by
\begin{equation}\label{f:d}
\dist(x) = \delta_{\R^n\setminus\Omega}(x) =
\inf_{y\in\partial\Omega} \pgauge(x-y),
\qquad x\in\overline{\Omega}\,.
\end{equation}
\end{dhef}

For every $x\in\overline{\Omega}$
we shall denote by $\proj(x) = \proj_{\partial\Omega}(x)$
the set of projections of $x$ in $\partial\Omega$.
For every $x_0\in\partial\Omega$
we define the
\textsl{proximal normal cone}
to $\partial\Omega$ at $x_0$ as the set
\begin{equation}\label{f:pcone}
N^P(x_0) :=\{
\lambda\, v;\
x_0+v\in\Omega,\
x_0\in\proj(x_0+v),\
\lambda\geq 0\}
\end{equation}
(see \cite{CSW}).
In general, $N^P(x_0)$ is convex but it may not
be closed.

\begin{dhef}[Singular set]
We say that $x\in\Omega$ is a
\textsl{regular point} of $\Omega$
if $\proj(x)$ is a singleton.
We say that $x\in\Omega$ is a
\textsl{singular point} of $\Omega$
if $x\in\Omega$ is not a regular point.
We denote by $\Sigma\subseteq\Omega$ the set
of all singular points of $\Omega$.
\end{dhef}

{}From Proposition~\ref{p:basic}(i),
$\Sigma$ coincides with the set of points
in $\Omega$ at which $\dist$ is not
differentiable.

Since $\Omega$ does not contain half spaces,
it is clear that, for every $x\in\Omega$,
the quantity $\leno(x)$ introduced in (\ref{f:leno})
is finite, hence the point $\cut(x)$ (see (\ref{f:cut}))
is well defined.

Given $x\in\Omega$, we call a \textsl{ray through $x$}
any closed segment $[y, \cut(x)]$ where $y\in\proj(x)$.
Thus there is a unique ray through $x\in\Omega$ if and only if
$\proj(x)$ is a singleton, that is,
if and only if $x\not\in\Sigma$.
From Proposition~\ref{p:raylength}
it is clear that, in any case, all rays through $x\in\Omega$
have the same length $\leno(x)$ in the
Minkowski distance.
Observe that
the length of a ray $[y,\cut(x)]$ is measured
by $\pgauge(\cut(x)-y)$ and,
in general, it is different from the length
$\pgauge(y-\cut(x))$ of $[\cut(x),y]$.

In Lemma~\ref{l:lusc} below we prove that
$\leno$ is an upper semicontinuous function in $\Omega$.
Furthermore, in Proposition~\ref{p:lusc}
we shall define an upper semicontinuous
extension $\len$ of $\leno$ to
$\overline{\Omega}$.

\begin{lemma}\label{l:lusc}
The function $\leno\colon\Omega\to\R$ defined in
(\ref{f:leno}) is upper semicontinuous in $\Omega$.
\end{lemma}

\begin{proof}
We follow the lines of the proof of
Proposition~{3.7} in \cite{EH}, where the same property
corresponding to the Euclidean distance is proved.
By contradiction,
assume that there exists a sequence $(x_k)_k\subset\Omega$
converging to a point $x\in\Omega$, such that
\begin{equation}\label{f:xc}
\lim_k \leno(x_k) > \leno(x)\,.
\end{equation}
For every $k\in\N$, choose $y_k\in\proj(x_k)$
and define
\[
\lambda_k = \min\left\{\frac{\leno(x_k)}{\dist(x_k)}\,,\,
2\frac{\leno(x)}{\dist(x)}\right\}\,.
\]
The sequence $(y_k)_k\subset\partial\Omega$
is clearly bounded.
Moreover,
since $\leno(x_k)\geq\dist(x_k)$ for every $k\in\N$
(see Remark~\ref{r:leno}), we have
$1\leq \lambda_k\leq 2\,\leno(x)/\dist(x)$,
hence the sequence $(\lambda_k)_k$ is also bounded.
We can then extract a subsequence, which we do not relabel,
such that
$\lim_k\lambda_k = \lambda\geq 1$,
$\lim_k y_k = y\in \partial \Omega$.
Since, by definition of $\leno(x_k)$,
$\lambda_k\leq \leno(x_k)/\dist(x_k)$,
we have that
$y_k\in\proj(y_k + \lambda_k (x_k-y_k))$ for every $k\in\N$.
Hence, by Lemma~\ref{l:pusc}, we infer that
$y\in\proj(y+\lambda(x-y))$.
By the definition of $\leno(x)$ we conclude that
$\lambda\leq \leno(x)/\dist(x)$,
so that $\lambda_k = \leno(x_k)/\dist(x_k)$
for $k$ large enough.
Then
\[
\lambda = \lim_k \frac{\leno(x_k)}{\dist(x_k)}\,,
\]
that is
\[
\lim_k\leno(x_k) = \lambda\, \lim_k\dist(x_k)
\leq \leno(x),
\]
in contradiction with (\ref{f:xc}).
\end{proof}

\begin{prop}\label{p:lusc}
The function $\len\colon\overline{\Omega}\to\R$ defined by
\begin{equation}\label{f:len}
\len(x) =
\begin{cases}
\leno(x)
&\textrm{if $x\in\Omega$},\\
\sup\{\dist(z);\ z\in\overline{\Omega}\ \textrm{and}\ x\in\proj(z)\}
&\textrm{if $x\in\partial\Omega$},
\end{cases}
\end{equation}
is upper semicontinuous in $\overline{\Omega}$.
Here $\leno\colon\Omega\to (0,+\infty)$ is the
function defined in (\ref{f:leno}).
\end{prop}

\begin{proof}
Recalling Lemma~\ref{l:lusc}
we have only to prove that $\len$ is upper
semicontinuous at every point of $\partial\Omega$.
Let $x_0\in\partial\Omega$.
Assume by contradiction that there exists
a sequence $(x_k)_k\subset\overline{\Omega}$,
converging to $x_0$,
such that
\begin{equation}\label{f:xxc}
\lim_k \len(x_k) > \len(x_0)\,.
\end{equation}
Upon passing to a subsequence,
it is enough to consider the following two cases:
\begin{itemize}
\item[(a)]
$x_k\in\partial\Omega$ for every $k\in\N$;
\item[(b)]
$x_k\in\Omega$ for every $k\in\N$.
\end{itemize}
In case (a), by the very definition of $\len$,
for every $k\in\N$ there exists a point $z_k$
such that
\begin{equation}\label{f:zk}
z_k\in\overline{\Omega},\quad
x_k\in\proj(z_k),\quad
\dist(z_k) > \len(x_k) - 1/k
\qquad (k\geq 1)\,.
\end{equation}
Since $(z_k)_k$ is a bounded sequence,
there exists a subsequence, which we do not relabel,
and a point $z_0\in\overline{\Omega}$ such that
$\lim_k z_k = z_0$.
From the upper semicontinuity of $\proj$
we have that $x_0\in\proj(z_0)$, hence
by (\ref{f:zk})
\[
\len(x_0)\geq\dist(z_0)\geq
\limsup_k\len(x_k)\,,
\]
in contradiction with (\ref{f:xxc}).

\noindent
Consider now case (b).
For every $k\in\N$, let $y_k\in\proj(x_k)$.
Up to a subsequence we can assume that
$(y_k)_k$ is convergent.
By the upper semicontinuity of $\proj$,
$(y_k)_k$ converges to $x_0$,
since $\proj(x_0) = \{x_0\}$.
Let us define
\[
\xi_k = \frac{x_k-y_k}{\pgauge(x_k-y_k)}\,,
\qquad k\in\N.
\]
By the characterization (\ref{f:leno2}) of $\len$,
for every $k$ we can choose $\lambda_k > \len(x_k) - 1/k$,
$\lambda_k>0$,
such that the point $z_k = y_k+\lambda_k \,\xi_k$
satisfies $y_k\in\proj(z_k)$.
Since $(\lambda_k)_k$ is bounded and
$\xi_k\in\partial K^0$ for every $k$,
we can extract another subsequence
(which we do not relabel)
such that $\lim_k \xi_k = \xi\in\partial K^0$
and $\lim_k\lambda_k = \lambda$.
Then we have $\lim_k z_k = x_0+\lambda\,\xi$ and,
by the upper semicontinuity of $\proj$,
$x_0\in\proj(x_0+\lambda\,\xi)$,
hence $\len(x_0)\geq \dist(x_0+\lambda \xi)$.
Collecting all the information we obtain
\[
\len(x_0)\geq\dist(x_0+\lambda\,\xi) = \lambda
=\lim_k\lambda_k\geq\limsup_k\len(x_k),
\]
which contradicts (\ref{f:xxc}).
\end{proof}

\begin{dhef}
The point $\cut(x)$, defined in (\ref{f:cut}),
is called the \textsl{ridge point}
of $x\in\Omega$.
The set
\[
\ridge = \{x\in\overline{\Omega};\ \dist(x) = \len(x)\}
\]
is called the \textsl{ridge}
of $\Omega$.
\end{dhef}

Observe that, if $x\in\ridge$,
then either $x=\cut(x)$,
or $x\in\partial\Omega$ and
$x\not\in\proj(y)$ for every $y\in\Omega$.
If $x\in\Omega$ is not a ridge point, then $\len(x)>\dist(x)$,
and,
from Proposition~\ref{p:basic}(iii),
the set $\proj(x)$ is a singleton,
hence $\dist$ is differentiable at $x$.
Thus the singular set $\Sigma$ is
contained in the ridge $\ridge$ of $\Omega$.
We remark that, in general,
the two sets do not coincide,
as it is shown in the following example.

\begin{example}\label{e:elli}
Let $\Omega = \{(x,y\in\R^2;\ x^2/a^2+y^2/b^2 < 1\}$,
where $0< b < a$ and $K=\overline{B}_1(0)$.
The points $P = ((a^2-b^2)/a, 0)$ and
$Q = (-(a^2-b^2)/a, 0)$ are the centers of curvature
of $\partial\Omega$ at $(a,0)$ and $(-a,0)$ respectively.
It can be checked that
$\ridge = [Q, P]$ whereas $\Sigma = (Q, P)$.
\end{example}

\section{Distance from the boundary of a smooth set}

Throughout the rest of the paper,
we assume that
\begin{equation}\label{f:Omega}
\Omega\subset\R^n\
\textrm{ is a nonempty, bounded, open connected set
of class $C^2$}.
\end{equation}

For every $x_0\in\partial\Omega$ we denote
by $\curv_1(x_0),\ldots, \curv_{n-1}(x_0)$ the principal
curvatures of $\partial\Omega$ at $x_0$,
and by $\nor(x_0)$ the \textsl{inward} normal unit vector
to $\partial\Omega$ at $x_0$.

Since the boundary of $\Omega$ is regular,
we will be able to extend $\dist$ outside $\overline\Omega$
in such a way that this extension turns out to be of
class $C^2$ in a tubular neighborhood of $\partial\Omega$
(see Theorem~\ref{t:regd} below).
This fact will allow us to define $D\dist$ and
$D^2\dist$ on points of $\partial\Omega$.
The extension of $\dist$ to $\R^n$ can be constructed as follows.
Let us define
\begin{equation}\label{f:distminus}
\dist^-(x) = \inf_{y\in\partial\Omega} \pgauge(y-x),
\qquad x\in\R^n\setminus\Omega\,,
\end{equation}
and the \textsl{signed distance function} from $\partial\Omega$
\begin{equation}\label{f:signeddist}
\sdist(x) =
\begin{cases}
\dist(x),
&\textrm{if $x\in\overline{\Omega}$},\\
-\dist^-(x),
&\textrm{if $x\in\R^n\setminus\overline{\Omega}$}\,.
\end{cases}
\end{equation}
For every $x\in\R^n\setminus{\Omega}$
define the set of projections from outside $\Omega$
\[
\proj^-(x) = \{y\in\partial\Omega;\
\dist^-(x) = \pgauge(y-x)\}\,,
\]
and we extend the projection operator $\Pi$ to be
\[
\proj^s(x) =
\begin{cases}
\proj(x),
&\textrm{if $x\in\overline{\Omega}$},\\
\proj^-(x),
&\textrm{if $x\in\R^n\setminus\overline{\Omega}$}\,.
\end{cases}
\]
We define also the set $\Sigma^-$ of singular points
of $\dist^-$ in $\R^n\setminus\overline{\Omega}$
and the set $\Sigma^s$ of singular points of
$\sdist$ in $\R^n$.
We clearly have the inclusion
$\Sigma^s\subseteq\Sigma\cup\Sigma^-\cup\partial\Omega$.

It is worth to observe that, in general,
$\proj^-(x)\neq \proj_{\overline{\Omega}}(x)$
for $x\in\R^n\setminus\overline{\Omega}$,
since $\pgauge$ need not be symmetric.

It is clear from the definition that
$\dist^-$ is the Minkowski distance from $\partial\Omega$
induced by the gauge function of $-K^0$.
We can interpret $\dist^-(x)$ as the Minkowski distance
of a point $x\in\R^n\setminus\Omega$ \textsl{to} the
boundary of $\Omega$.

As in (\ref{f:len}),
we define a function $\len^-$ on $\partial\Omega$ by
\begin{equation}\label{f:lenm}
\len^-(x_0) =
\sup\{\dist^-(z);\ z\in\R^n\setminus{\Omega}\ \textrm{and}\ x_0\in\proj^-(z)\}
\qquad x_0\in\partial\Omega\,.
\end{equation}
In the following lemma we show that,
under our regularity assumption on $\Omega$,
the functions $\len$ and $\len^-$ are bounded from below
by a positive constant on $\partial\Omega$.

\begin{lemma}\label{l:taupos2}
Let $\Omega$ satisfy $(\ref{f:Omega})$.
Then there exists a positive constant $\mu$ such that
$\len(x_0)$, $\len^-(x_0) \geq \mu$ for every $x_0\in\partial\Omega$.
As a consequence, $\overline{\ridge}\subseteq \Omega$.
\end{lemma}

\begin{proof}
We shall only prove the assertion concerning $\len$,
the other being similar.
Since $\Omega$ is of class $C^2$,
then it satisfies
a uniform interior sphere condition of radius $r>0$,
that is,
for every boundary point $x_0\in\partial\Omega$
there exists $z_0\in\Omega$ such that
$x_0\in \overline{B}_r(z_0)\subset\overline{\Omega}$.
Let $R>0$ denote
the maximum of the principal radii of curvature
of $\partial K^0$,
so that
$K^0$ (and $-K^0$) slides freely inside $\overline{B}_R$,
that is,
for every boundary point $y\in\partial B_R$
there exists $z\in B_R$ such
that $y\in z- K^0\subset\overline{B}_R$
(see \cite[Corollary~3.2.10]{Sch}).
We are going to prove that, for $\mu = r/R$,
$-\mu\, K^0$ slides freely inside $\Omega$.
Let $x_0\in\partial\Omega$.
We have to show that
there exists $z\in\Omega$ such that
$x_0\in z- \mu\,K^0
\subset\overline{\Omega}$.
From the uniform interior sphere condition,
there exists $z_0\in\Omega$
such that
$x_0\in \overline{B}_r(z_0)\subset\overline{\Omega}$.
On the other hand, $-\mu\, K^0$
slides freely inside
$\mu\, \overline{B}_R = \overline{B}_r$.
Thus, there exists $z\in\overline{B}_r(z_0)$ such that
$x_0\in z-\mu\, K^0\subset
\overline{B}_r(z_0)$.
We have proved that for every $x_0\in\partial \Omega$
there exists $z\in \Omega$ such that $x_0\in z-\mu K^0\subset \overline\Omega$.
Then $\mu=\pgauge(z-x_0)\leq \pgauge(z-y)$ for every $y\in \partial \Omega$,
that is $x_0\in\proj(z)$ and $\dist(z)=\mu$.
Hence $\len(x_0) \geq \mu$.
\end{proof}

\begin{remark}\label{r:ridins}
As a consequence of Lemma~\ref{l:taupos2}, we obtain
that, for $x_0\in\partial \Omega$, the supremum in the definition (\ref{f:len})
of $\len(x_0)$ is achieved by a point $z\in\Omega$.
Namely, from the continuity of $\dist$ and the upper
semicontinuity of $\proj$, the supremum is achieved
at a point $z\in\overline{\Omega}$.
Should $z\in\partial\Omega$, then $x_0\in\proj(z)$
would imply $z=x_0$ and $\len(x_0) = 0$,
in contradiction with Lemma~\ref{l:taupos2}.

In the definition (\ref{f:lenm}) of $\len^-(x_0)$ the
supremum need not be achieved, but the argument above
shows that we can replace $\R^n\setminus\Omega$
by $\R^n\setminus\overline{\Omega}$ in the definition
of $\len^-$.
\end{remark}

The following result
relates the gradient of $\sdist$ to the
inward normal of $\partial\Omega$.

\begin{lemma}\label{l:pis}
Let
$x\in\R^n\setminus\Sigma^s$,
$x\not\in\partial\Omega$,
and let $\proj^s(x) = \{x_0\}$.
Then
\begin{equation}\label{f:diffd}
D\sdist(x) = \frac{\nor(x_0)}{\gauge(\nor(x_0))}\,.
\end{equation}
\end{lemma}

\begin{proof}
Consider first the case $x\in{\Omega}\setminus{\Sigma}$.
{}From Lemma~3.5 in \cite{Pi} we have that
$D\dist(x) = \lambda\,\nor(x_0)$ for some $\lambda \geq 0$.
On the other hand, from Proposition~\ref{p:basic}(i)
we have that $D\dist(x) = D\pgauge(x-x_0)$, and,
from Theorem~\ref{t:sch}(iii), this vector
belongs to $\partial K$.
Then
$1 = \gauge(D\dist(x)) = \lambda\, \gauge(\nor(x_0))$,
that is $\lambda = 1/ \gauge(\nor(x_0))$
and (\ref{f:diffd}) follows.

\noindent
Consider now the case
$x\in\R^n\setminus(\overline{\Omega}\cup\Sigma^-)$.
{}From Proposition~\ref{p:basic}(i)
and Lemma~3.5 in \cite{Pi}
we have that
$D\gau{(-K)^0}(x-x_0) = D\dist^-(x) = -\lambda\,\nor(x_0)$,
for some $\lambda\geq 0$.
Upon observing that
\[
(-K)^0 = -K^0,\qquad
\gau{-K}(\xi) = \gau{K}(-\xi),\quad\forall\xi\in\R^n,
\]
we deduce that $D\sdist(x) = D\pgauge(x_0-x)\in\partial K$,
hence
\[
1=
\gauge(D\sdist(x))
= \lambda\, \gauge(\nor(x_0))\,,
\]
and (\ref{f:diffd}) follows.
\end{proof}

We are now in a position to characterize the normal
directions to the boundary of $\Omega$.

\begin{prop}\label{p:pnc}
For every $x_0\in\partial\Omega$,
the proximal normal cone of $\Omega$ at $x_0$,
defined in (\ref{f:pcone}),
is given by
\begin{equation}\label{f:pnc}
N^P(x_0) = \{\lambda\, D\gauge(\nor(x_0));\ \lambda\geq 0\}\,,
\end{equation}
whereas the proximal normal cone of
$\R^n\setminus\overline{\Omega}$ at $x_0$
with respect to $\dist^-$
is $-N^P(x_0)$.
\end{prop}

\begin{proof}
By Proposition~\ref{p:basic}(i) and Lemma~\ref{l:pis}
the vector $D\gauge(\nor(x_0))$ belongs to $N^P(x_0)$.
On the other hand,
by the definition of $N^P(x_0)$,
we have that $w\in N^P(x_0)\setminus\{0\}$
if and only if there exists $\mu > 0$ such that
$x_0\in\proj(x_0+\mu\, w)$.
Then, for $\epsilon\in (0,\mu)$,
from Proposition~\ref{p:basic}(iii)
we have $\proj(x_0+\epsilon\, w) = \{x_0\}$.
From Lemma~\ref{l:pis}
and Proposition~\ref{p:basic}(i)
we have that
\[
D\gauge(\nor(x_0)) =
D\gauge(D\dist(x_0+\epsilon\, w)) = \frac{w}{\pgauge(w)}\,,
\]
hence $w = \lambda\, D\gauge(\nor(x_0))$
with $\lambda = \pgauge(w) > 0$.
The computation of
the proximal normal cone of
$\R^n\setminus\overline{\Omega}$ at $x_0$
is similar.
\end{proof}

\begin{remark}\label{r:pnc}
As a consequence of Proposition~\ref{p:pnc} and Proposition~\ref{p:basic}(i)
it is clear that,
{}from any point $x_0\in\partial\Omega$,
$D\gauge(\nor(x_0))$ is the unique inward ``normal'' direction
with the properties
\[
\proj^s[x_0+t\,D\gauge(\nor(x_0))] = \{x_0\},\qquad
\sdist(x_0+t\,D\gauge(\nor(x_0))) = t
\]
for $t$ small enough
(see also \cite[Lemma~2.2]{LN}).
In the following proposition
we shall show that, under our assumptions,
the relations above hold for
$t\in(-\len^-(x_0),\len(x_0))$.
This is a well known fact in Riemannian geometry
(see \cite[\S III.4]{Sa}).
\end{remark}

\begin{prop}\label{p:conj}
For every $x_0\in\partial\Omega$,
we have that
\begin{equation}\label{f:conj}
\proj^s(x_0+t\, D\gauge(\nor(x_0))) = \{x_0\}
\qquad
\forall t\in (-\len^-(x_0), \len(x_0)),
\end{equation}
whereas $x_0\not\in\proj^s(x_0+t\, D\gauge(\nor(x_0)))$
for every $t<-\len^-(x_0)$
and for every $t>\len(x_0)$.
\end{prop}

\begin{proof}
We prove the assertion for $t\geq 0$, the case $t\leq 0$
being similar.
Fixed $x_0\in\partial \Omega$, by Remark~\ref{r:ridins}
there exists $z\in \Omega$ such that $x_0\in\proj(z)$
and $\dist(x_0)=\len(x_0)$.
By Proposition~\ref{p:basic}(iii),
$\{x_0\}=\proj(x)$ for every $x$ in the
segment $(x_0,z)$.
Moreover, by Proposition~\ref{p:basic}(i)
and (\ref{f:diffd}),
the segment $(x_0,z)$ can be parameterized by
\[
x_t = x_0+t\, \frac{z-x_0}{\pgauge(z-x_0)} =
x_0 + t\, D\gauge(\nor(x_0)),\quad
t\in (0,\len(x_0))\,,
\]
so that (\ref{f:conj}) holds true.
On the other hand, by the very definition
of $\len(x_0)$,
it cannot happen that
$x_0\in\proj^s(x_t)$ and $t=\dist(x_t)>\len(x_0)$.
\end{proof}

The following result, which will be used in the sequel,
states that every point $x\in\Omega\setminus\overline{\Sigma}$
belongs to the interior of the ray
$(p(x), \cut(x))$, where $p(x)$ is the unique projection of
$x$ on $\partial\Omega$,
and $\cut(x)$ is the ridge point of $x$
defined in (\ref{f:cut}).

\begin{lemma}\label{l:Borsuk}
Let $x\in\Omega\setminus\overline{\Sigma}$ and let
$p(x)$ denote its unique projection on $\partial\Omega$.
Then there exists a point $x_1\in\Omega\setminus\overline{\Sigma}$
such that $\proj(x_1) = \{p(x)\}$ and
$x\in (p(x), x_1)$.
\end{lemma}

\begin{proof}
Since $\Omega\setminus\overline{\Sigma}$ is an open set,
there exists $r>0$ such that
$B_r(x)\subset \Omega\setminus\overline{\Sigma}$.
For every $z\in B_r(x)$ let $p(z)$ denote the unique
projection of $z$ on $\partial\Omega$.
{}From Lemma~\ref{l:pusc}, the map $p$ is continuous
in $B_r(x)$.
Let $V\subset\partial\Omega$ be a local chart on
$\partial\Omega$ containing $p(x)$, and let
$\varphi\colon V\to \R^{n-1}$ be a local coordinate system,
that is, a $C^2$ bijection from $V$ to the open
set $\U = \varphi(V)\subset\R^{n-1}$.
Since $p$ is continuous in $B_r(x)$, there exists
$\delta\in (0,r)$ such that
$p(\overline{B}_{\delta}(x))\subset V$.
Consider now the map
$\psi\colon S^{n-1}\to\R^{n-1}$ defined by
\[
\psi(u) = \varphi(p(x+\delta\, u))\,,
\qquad u\in S^{n-1}\,.
\]
It is clear that this map is continuous,
being the composition of continuous maps.
{}From the Borsuk-Ulam theorem
(see \cite[Corollary~4.2]{Dei})
there exists $u_0\in S^{n-1}$ such that
$\psi(u_0) = \psi(-u_0)$.
Since $\varphi$ is one-to-one,
we deduce that $p(x+\delta\, u_0) = p(x-\delta\, u_0)$,
that is, the points $x_1 = x+\delta\, u_0$
and $x_2 = x-\delta\, u_0$ have the same projection
$x_0\in\partial\Omega$.
{}From Proposition~\ref{p:basic}(i)
and Lemma~\ref{l:pis} we have that
\[
x_i = x_0+\dist(x_i)D\gauge(\nor(x_0))\,,\qquad
i=1,2.
\]
It is clear that the points
$x_0$, $x_1$ and $x_2$ are collinear.
Assume, just to fix the ideas, that
$\dist(x_1) > \dist(x_2)$, that is,
$\pgauge(x_1-x_0) > \pgauge(x_2-x_0)$,
in such a way that $x_2\in (x_0, x_1)$.
Since $x\in (x_1, x_2)$
by definition of $x_1$ and $x_2$,
we also have $x\in (x_0, x_1)$.
Finally, recalling that $\proj(x_1) = \{x_0\}$, from
Proposition~\ref{p:basic}(iii) we conclude that
$p(x) = x_0$.
\end{proof}

\begin{remark}\label{r:bors}
As consequence of Lemma \ref{l:Borsuk},
we have that $\len(x)>\dist(x)$ for every
$x\in \Omega \setminus \overline{\Sigma}$.
Hence
$\ridge \subseteq \overline{\Sigma}$.
\end{remark}

For many calculations we shall use
a preferred system of coordinates
in order to parameterize $\partial\Omega$
in a neighborhood of a point $x_0$.

\begin{dhef}[Principal coordinate system]\label{d:coord}
Let $x_0\in\partial\Omega$.
We call \textsl{principal coordinate system} at $x_0$
the coordinate system such that $x_0 = 0$, $e_n = \nor(x_0)$
and $e_i$ coincides with the $i$-th principal direction
of $\partial\Omega$ at $x_0$, $i=1,\ldots,n-1$.
\end{dhef}

Using the principal coordinate system at $x_0$,
$\partial\Omega$ can be parameterized in a neighborhood of
$x_0$ by a map
\[
X\colon\U\to\R^n,\quad
X(y) = (y,\phi(y)),
\]
where $\U\subset\R^{n-1}$ is a neighborhood of the origin,
and $\phi\colon\U\to\R$ is a map of class $C^2$
satisfying
\begin{equation}\label{f:fiz}
\phi(0) = 0,\quad
\frac{\partial\phi}{\partial y_i}(0) = 0,\quad
i=1,\ldots,n-1.
\end{equation}
We shall refer to $X$ as a
\textsl{standard parametrization} of
$\partial\Omega$ in a neighborhood of $x_0$.
If we denote by $N(y) = \nor(X(y))$, $y\in\U$,
then
\[
N(y) = \frac{1}{\sqrt{1+|D_y\phi(y)|^2}}\,
(-D_y\phi(y), 1),
\qquad y\in\U,
\]
where
$D_y\phi = (\partial\phi/\partial y_1,\ldots,
\partial\phi/\partial y_{n-1})$.
Since, in the principal coordinate system,
the vectors $e_1,\ldots,e_{n-1}$ are the
principal directions of $\partial\Omega$ at $x_0$, we have
\begin{equation}\label{f:N}
N(0) = e_n,\quad
\frac{\partial N}{\partial y_i}(0) = -\curv_i e_i,\quad
i=1,\ldots,n-1,
\end{equation}
where $\curv_1,\ldots,\curv_{n-1}$ are the principal curvatures
of $\partial\Omega$ at $x_0$.
Starting from the identities
\[
\frac{\partial \phi}{\partial y_j}=-N_j(y) \sqrt{1+|D_y\phi(y)|^2},\quad j=1,\ldots,n-1\,,
\]
and differentiating, we easily obtain
\begin{equation}\label{f:fiz2}
\frac{\partial^2 \phi}{\partial y_i\partial y_j}(0) =
-\frac{\partial N_j}{\partial y_i}(0)
=\curv_i\, \delta_{ij},\qquad i, j=1,\ldots,n-1,
\end{equation}
where $\delta_{ij}$ denotes
the Kroneker symbol.
Moreover, since
\[
\frac{\partial X}{\partial y_i}(y)=e_i+\left(0,\frac{\partial\phi}{\partial y_i}(y)\right),
\qquad
\frac{\partial^2 X}{\partial y_i\partial y_j}(y)=
\left(0,\frac{\partial^2\phi}{\partial y_i\partial y_j}(y)\right),
\]
for every $i, j=1,\ldots,n-1$, by (\ref{f:fiz}) and (\ref{f:fiz2}), we get
\begin{equation}\label{f:X}
\frac{\partial X}{\partial y_i}(0) = e_i,\
\frac{\partial^2 X}{\partial y_i\partial y_j}(0) =
\curv_i \delta_{ij}\,e_n,\quad
i,j=1,\ldots,n-1\,.
\end{equation}

In the following lemma we introduce the main tool of
our theory, that is a parametrization of $\Omega$ that
will allow us to prove the regularity of the signed distance
near $\partial \Omega$, as well as the regularity of the set
$\overline{\Sigma}$.
Moreover, in Section~\ref{s:PDE} we shall
use this parametrization in order to prove a
change of variables formula in multiple integrals on $\Omega$.

In the following computations, we extend the notion of inward normal
for $x$ in a tubular neighborhood of $\partial \Omega$,
by setting $\nor(x)$ to be the gradient of the Euclidean signed distance
of $x$ from $\partial \Omega$. In this way $D\nor(x_0)$ is well
defined for every $x_0\in \partial\Omega$, and $D\nor(x_0)\nor(x_0)=0$.

\begin{lemma}\label{l:detpsi}
Let $Y\colon\U\to\R^n$, $\U\subset\R^{n-1}$ open,
be a local parametrization of $\partial\Omega$ of class $C^2$.
Let $\Psi\colon \U\times\R\to\R^n$ be the map
defined by
\begin{equation}\label{f:psi}
\Psi(y,t) = Y(y) + t D\gauge(\nor(Y(y))),
\quad (y,t)\in \U\times\R\,.
\end{equation}
Then $\Psi\in C^1(\U\times\R)$, and
\begin{equation}\label{f:calcdetpsi}
\det D\Psi(y,t) = \gauge(\nor(Y(y))\,
\sqrt{g(y)}\,
\det[I+t\, D^2\gauge(\nor(Y(y)))\, D\nor(Y(y))]
\end{equation}
for every
$(y,t)\in\U\times\R$,
where
\[
g(y) = \det (g_{ij}(y)),
\qquad
g_{ij}(y) = \pscal{\frac{\partial Y}{\partial y_i}(y)}%
{\frac{\partial Y}{\partial y_j}(y)}
\,,\quad
i,j=1,\ldots,n-1,
\]
is the determinant of the matrix of the metric coefficients.
\end{lemma}

\begin{proof}
Since $D\gauge\in C^1(\R^n\setminus\{0\})$
and $\nor\circ Y\in C^1(\U)$,
$\Psi$ is of class $C^1$ in $\U\times\R$.
Let us fix $y_0\in \U$, let $x_0 = Y(y_0)$,
and define
\[
Q = I + t\, D^2\gauge(\nor(x_0))\, D\nor(x_0),
\quad
w_i = \frac{\partial Y}{\partial y_i}(y_0),
\quad i = 1,\ldots,n-1.
\]
We have that
\[
\det D\Psi(y_0,t) =
\det\left[Q\,w_1,\ldots,
Q\,w_{n-1}\,, D\gauge(\nor(x_0))
\right]\,.
\]
Notice that,
for every $i=1,\ldots,n-1$,
$\pscal{w_i}{\nor(x_0)} = 0$ and,
thanks to (\ref{f:eigenzero}),
\[
\pscal{D^2\gauge(\nor(x_0))\, D\nor(x_0)\, w_i}{\nor(x_0)} =
\pscal{D\nor(x_0)\, w_i}{D^2\gauge(\nor(x_0))\,\nor(x_0)} = 0\,.
\]
Hence $\pscal{Q\, w_i}{\nor(x_0)} = 0$
for every $i=1,\ldots,n-1$,
and $Q\,\nor(x_0) = \nor(x_0)$.
Upon observing that
$D\gauge(\nor(x_0)) = w
+ \gauge(\nor(x_0))\, \nor(x_0)$,
with $\pscal{w}{\nor(x_0)} = 0$, we get
\[
\begin{split}
\det D\Psi(y_0,t) & =
\gauge(\nor(x_0))\,\det\left[Q\,w_1,\ldots,
Q\,w_{n-1}\,, \nor(x_0)
\right]
\\ & =
\gauge(\nor(x_0))\,\det Q\,
\det\left[w_1,\ldots,w_{n-1}, \nor(x_0)
\right]\,.
\end{split}
\]
Finally, the matrix
$C = \left[w_1,\ldots,w_{n-1}, \nor(x_0)\right]$
satisfies
\[
\begin{split}
(C^T\, C)_{ij} & = \pscal{w_i}{w_j},
\quad i,j=1,\ldots,n-1,\\
(C^T\, C)_{kn} & = (C^T\, C)_{nk} = \delta_{kn},
\quad k=1,\ldots,n\,.
\end{split}
\]
Hence
$(\det C)^2 = \det(C^T\, C) = g(y_0)$,
and (\ref{f:calcdetpsi}) is proved.
\end{proof}

The crucial point in what follows is to have
exact information about the degeneracy of the
map $\Psi$. The first step in this direction
is to simplify the computation of
$\det[I+t\, D^2\gauge(\nor)\, D\nor]$.

\begin{lemma}\label{l:detpsi2}
For every $x_0\in\partial\Omega$,
we have that
\begin{equation}\label{f:calcdetpsi2}
\det[I+t\, D^2\gauge(\nor(x_0))\, D\nor(x_0)]
= \det(I_{n-1} - t R\, D),
\end{equation}
where
$D$, $R\in\matr_{n-1}$ are the symmetric
square matrices with entries
\begin{equation}\label{f:DR}
D_{ij} = -\curv_i(x_0)\,\delta_{ij}\,,
\quad
R_{ij} = \pscal{D^2\gauge(\nor(x_0))\,e_j}{e_i},
\quad
i,j = 1,\ldots,n-1,
\end{equation}
in the principal coordinate system at $x_0$.
\end{lemma}

\begin{proof}
It is enough to observe that,
in the principal coordinate system at $x_0$,
$\pscal{D\nor(x_0)\,e_i}{e_j} = -\curv_i\,\delta_{ij}$
for every $i,j=1,\ldots,n-1$,
whereas, from (\ref{f:eigenzero}),
\[
\pscal{D^2\gauge(\nor(x_0))\, D\nor(x_0)\,e_i}{\nor(x_0)}
= \pscal{D\nor(x_0)\,e_i}{D^2\gauge(\nor(x_0))\, \nor(x_0)}
= 0
\]
for every $i=1,\ldots,n-1$,
hence (\ref{f:calcdetpsi2}) follows.
\end{proof}

The equality (\ref{f:calcdetpsi2}) is fruitful,
since it allows us to deal with a positive definite matrix
$R$ with known inverse matrix, as stated
in the following result.

\begin{lemma}\label{l:detpsi3}
Let $x_0\in\partial\Omega$,
and let $H\in\matr_{n-1}$ be the matrix
with entries
\begin{equation}\label{f:H}
H_{ij} = \pscal{D^2\pgauge(D\gauge(\nor(x_0)))\,e_j}{e_i},
\qquad
i,j = 1,\ldots,n-1,
\end{equation}
in the principal coordinate system at $x_0$.
Then
\begin{equation}\label{f:RH}
\gauge(\nor(x_0))\, R\, H = I_{n-1}\,,
\end{equation}
where $R\in\matr_{n-1}$ is the matrix
defined in (\ref{f:DR}).
As a consequence, $R$ and $H$ are positive
definite matrices.
\end{lemma}

\begin{proof}
In order to prove (\ref{f:RH}),
let us differentiate the first
identity in (\ref{f:polare}) with respect to $\xi$,
obtaining
\[
D^2\pgauge(D\gauge(\xi))\, D^2\gauge(\xi) =
\frac{1}{\gauge(\xi)}\, I -
\frac{1}{\gauge(\xi)^2}\, \xi\otimes D\gauge(\xi),
\qquad\forall\xi\neq 0\,.
\]
{}From the symmetry of $D^2\pgauge$ and $D^2\gauge$,
the adjoint matrix identity can be written as
\[
D^2\gauge(\xi)\, D^2\pgauge(D\gauge(\xi))  =
\frac{1}{\gauge(\xi)}\, I -
\frac{1}{\gauge(\xi)^2}\, D\gauge(\xi)\otimes\xi,
\qquad\forall\xi\neq 0\,.
\]
Notice that, if $v$ is orthogonal to
$\nor = \nor(x_0)$, then
$(D\gauge(\nor)\otimes\nor)\,v = 0$, hence,
in the principal coordinate system,
for every $i,j=1,\ldots,n-1$,
\[
\begin{split}
\frac{1}{\gauge(\nor)}\,\delta_{ij}
& = \pscal{D^2\gauge(\nor)\, D^2\pgauge(D\gauge(\nor))\, e_i}{e_j}
\\ & =
\sum_{l=1}^{n-1} \pscal{D^2\gauge(\nor)\, e_l}{e_j}
\pscal{D^2\pgauge(D\gauge(\nor))\, e_i}{e_l}\,,
\end{split}
\]
where the last equality is due to the fact that,
from (\ref{f:eigenzero}),
$D^2\gauge(\nor)\, e_n =
D^2\gauge(\nor)\, \nor = 0$.
The relation above is exactly (\ref{f:RH}).

Recalling that $D^2\gauge(\nor)$ and $D^2\pgauge(D\gauge(\nor))$
are positive semidefinite (due to the convexity of
$\gauge$ and $\pgauge$),
it follows that $R$ and $H$ are both
positive semidefinite.
As a consequence of (\ref{f:RH}), $R$ and $H$ are also invertible,
hence they are positive definite.
\end{proof}

The fundamental step in the proof of the $C^2$ regularity
of $\sdist$ in a tubular neighborhood of $\partial\Omega$
and of $\dist$ in $\overline\Omega\setminus\overline{\Sigma}$
is the following.

\begin{theorem}\label{t:detpsi}
Let $Y\colon\U\to\R^n$, $\U\subset\R^{n-1}$ open,
be a local parametrization of $\partial\Omega$ of class $C^2$.
Let $\Psi\colon \U\times\R\to\R^n$ be the map
defined in (\ref{f:psi}).
Then $\det D\Psi(y,t) > 0$ for every
$y\in\U$ and every $t\in (-\len^-(Y(y)), \len(Y(y)))$.
\end{theorem}

\begin{proof}
Let us fix $x_0\in\partial \Omega$. There is no loss of
generality in assuming $x_0 = Y(0)$.
In the principal coordinate system at $x_0$,
we have also $x_0=0$.
We shall denote $\nor = \nor(x_0) = e_n$.
From Lemmas~\ref{l:detpsi}
and~\ref{l:detpsi2},
in this coordinate system we have to prove
that
\begin{equation}\label{f:detC}
\det(I_{n-1} - t R\, D) > 0\,,
\qquad\forall t\in (-\len^-(x_0), \len(x_0)),
\end{equation}
where
$D$, $R\in\matr_{n-1}$ are the
square matrices defined in (\ref{f:DR}).

Observe that, from (\ref{f:RH})
\[
I_{n-1} - t\, R\, D = \gauge(\nor)\, R\, H - t\, R\, D
={\gauge(\nor)}\,R\, \left(H-\frac{t}{\gauge(\nor)}\, D\right),
\]
where $H$ is the matrix defined
in (\ref{f:H}). Since $\det (R)>0$ by Lemma~\ref{l:detpsi3},
in order to prove (\ref{f:detC}),
it is enough to prove the following claim.

\smallskip\noindent
\textsl{Claim.}
The matrix $A(t)\in \matr_{n-1}$
defined by
\[
A(t) = H-\frac{t}{\gauge(\nor)}\, D,
\qquad t\in\R,
\]
is positive definite for every
$t\in (-\len^-(x_0),\len(x_0))$.

\smallskip\noindent
\textsl{Proof.}
For $t=0$ we have $A(0) = H$, which is positive definite
by Lem\-ma~\ref{l:detpsi3}.
We shall now prove the claim for $t\in (0,\len(x_0))$,
the case $t\in(-\len(x_0), 0)$ being similar.

\noindent
We start proving that,
for every $t\in (0,\len(x_0))$, $A(t)$ is
positive semidefinite.
Let $z_t = x_0 + t D\gauge(\nor)$.
From Proposition~\ref{p:conj}
we have that $\proj(z_t) = \{x_0\}$, and
$\dist(z_t) = t = \pgauge(z_t-x_0)$.
Moreover, the following ball in the Minkowski norm
\[
B = \{x\in\R^n;\ \pgauge(z_t-x) < t\}
\]
is contained in $\Omega$,
and $x_0\in\partial B \cap \partial\Omega$.
Thus
\[
h(y) := \pgauge(z_t - X(y)) \geq t
\quad
\forall y\in\U,
\quad h(0) = t\,.
\]
Then, the function $h\in C^2(\U)$ has a local minimum point
at $y=0$, hence its Hessian matrix at $y=0$ must be
positive semidefinite.
It is straightforward to check that
\begin{equation*}
\begin{split}
\frac{\partial h}{\partial y_i}(y)= &
\pscal{D\pgauge(z_t-X(y))}{\frac{\partial X}{\partial y_i}(y)}\,,\\
\frac{\partial^2 h}{\partial y_i\partial y_j}(y) = &
\pscal{D^2\pgauge(z_t-X(y))\frac{\partial X}{\partial y_j}(y)}
{\frac{\partial X}{\partial y_i}(y)} \\
 & +
\pscal{D\pgauge(z_t-X(y))}{\frac{\partial^2 X}{\partial y_i\partial y_j}(y)}\,.
\end{split}
\end{equation*}
Using (\ref{f:X}) we obtain
\[
\frac{\partial^2 h}{\partial y_i\partial y_j}(0) =
\pscal{D^2\pgauge(z_t-x_0) \, e_i}{e_j}
-\curv_i\,\delta_{ij}\,\pscal{D\pgauge(z_t-x_0)}{\nor}\,.
\]
Since $z_t-x_0 = t\, D\gauge(\nor)$, from the positive $0$-homogeneity of
$D\pgauge$, the positive $(-1)$-homogeneity of $D^2\pgauge$, and (\ref{f:polare})
we get
\[
D\pgauge(z_t-x_0) = D\pgauge(D\gauge(\nor))=\frac{\nor(y)}{\gauge(\nor(y))},\quad
D^2\pgauge(z_t-x_0) = \frac{1}{t}\, D^2\pgauge(D\gauge(\nor))\,,
\]
hence
\[
\frac{\partial^2 h}{\partial y_i\partial y_j}(0) =
\frac{1}{t}\, \pscal{D^2\pgauge(D\gauge(\nor)) \, e_i}{e_j}
-\frac{\curv_i}{\gauge(\nor)}\,\delta_{ij}\,,
\]
that is,
$D^2 h(0) = (1/t) H - (1/\gauge(\nor))\, D = A(t)/t$.
Since the Hessian of $h$ at $y=0$ is a positive semidefinite matrix,
we conclude that $A(t)$ is also a positive semidefinite matrix
for $t\in(0,\len(x_0))$.

\noindent
Let us prove that, in fact,
$A(t)$ is positive definite for every $t\in (0,\len(x_0))$.
It is plain that it is equivalent to prove that
$\lambda\,\gauge(\nor)\, H - D$
is positive definite for every $\lambda > 1/\len(x_0)$.
Let us fix $\lambda > 1/\len(x_0)$ and choose
$\lambda' \in (1/\len(x_0), \lambda)$,
so that $\lambda'\,\gauge(\nor)\, H - D$
is positive semidefinite.
Since,
by Lemma~\ref{l:detpsi3}, $H$ is positive definite,
\[
\begin{split}
\pscal{(\lambda\,\gauge(\nor)\, H - D)w}{w} & =
\pscal{(\lambda'\,\gauge(\nor)\, H - D)w}{w} +
(\lambda-\lambda')\pscal{\gauge(\nor)\, H\, w}{w}
\\ & \geq
(\lambda-\lambda')\pscal{\gauge(\nor)\, H\, w}{w}
> 0
\end{split}
\]
for every $w\neq 0$, that is,
$\lambda\,\gauge(\nor)\, H - D$
is positive definite.
This concludes the proof of the claim and of
the proposition.
\end{proof}

\begin{remark}
It is well known that $\dist$ is the unique viscosity solution of
the Hamilton-Jacobi equation $\gauge(D\dist) = 1$ in $\Omega$
vanishing on $\partial\Omega$
(see \cite{BaCD,CaSi,Li}).
The map $\Psi$ defined in (\ref{f:psi}) gives the characteristic curves
associated to this PDE.
The regularity of the map $\Psi$ proved in Theorem~\ref{t:detpsi}
will be used in
Theorem~\ref{t:regd} in order to prove the $C^2$ regularity of
$\sdist$ in a tubular neighborhood of $\partial\Omega$.
It should be noted that this kind of result in fact
follows from the local existence theory for first order PDEs
based on the method of characteristics.
Nevertheless
the regularity of the map $\Psi$
is an essential tool in order to prove the $C^2$ regularity of
$\dist$ in the whole set $\Omega\setminus\overline{\Sigma}$
(see Theorem~\ref{t:regd2} below).
\end{remark}

\begin{cor}\label{c:sigma}
The ridge set
$\ridge$ has zero Lebesgue
measure.
\end{cor}

\begin{proof}
Let $Y_k\colon\U_k\to\R^n$, $\U_k\subset\R^{n-1}$ open,
$k=1,\ldots,N$,
be local parameterizations of $\partial\Omega$ of class $C^2$,
such that $\bigcup_{k=1}^{N} Y_k(\U_k) = \partial\Omega$.
For every $k=1,\ldots,N$,
let $\Psi_k\colon \U_k\times\R\to\R^n$ be the map
\[
\Psi_k(y,t) = Y_k(y) + t D\gauge(\nor(Y_k(y))),
\quad (y,t)\in \U_k\times\R\,.
\]
For every $k=1,\ldots,n$ let $U_k\subset\U_k$
be a compact set such that $\bigcup_k Y_k(U_k)$ covers
$\partial\Omega$,
and let
\[
A_k = \{(y,t);\ y\in U_k,\ t\in [0,\len(Y_k(y))]\}\,.
\]
From Proposition~\ref{p:lusc},
$\len$ is an upper semicontinuous function,
hence for every $k=1,\ldots,n$,
$A_k$ is a compact set and
the Lebesgue measure of
the graph
\[
\Psi_k^{-1}(\ridge)\cap A_k =
\{(y,t)\in U_k\times\R;\ t=\len(Y_k(y))\}
\]
vanishes.
{}From Theorem~\ref{t:detpsi} we know that,
for every $k=1,\ldots,N$,
$\Psi_k\in C^1(\U_k\times\R)$,
hence it is Lipschitz continuous on
the compact set $A_k$.
Let $L$ be the maximum of the Lipschitz constants
of the functions $\Psi_1,\ldots,\Psi_N$.
Since $\bigcup_{k=1}^N \Psi_k(A_k) = \overline{\Omega}$,
and hence
$\ridge\subseteq \bigcup_{k=1}^N \Psi_k\left(\Psi_k^{-1}(\ridge)\cap A_k\right)$,
we finally get
\[
\meas(\ridge)\leq
\sum_{k=1}^N \meas\left[\Psi_k\left(
\Psi_k^{-1}(\ridge)\cap A_k\right)\right]
\leq L^n
\sum_{k=1}^N \meas\left[
\Psi_k^{-1}(\ridge)\cap A_k\right]
= 0
\]
and the proof is complete.
\end{proof}

\begin{theorem}[Regularity of $\sdist$]\label{t:regd}
The function $\sdist$ is of class $C^2$ in
a tubular neighborhood of $\partial\Omega$ of the form
\[
A_{\mu} = \{
x\in\R^n;\ -\mu < \sdist(x) < \mu\}
\]
for some $\mu > 0$.
Furthermore, for every $x_0\in\partial\Omega$
\begin{equation}\label{f:diffdb}
D\sdist(x_0+tD\gauge(\nor(x_0)))=D\sdist(x_0)\,,
\quad t\in(-\mu,\mu)\,,
\end{equation}
and then
the identity (\ref{f:diffd}) holds
also for $x\in\partial\Omega$.
If
$\Omega$ is of class $C^{k,\alpha}$
and $\gauge\in C^{k,\alpha}(\R^n\setminus\{0\})$,
for some $k\geq 2$ and $\alpha\in [0,1]$,
then $\sdist$ is of class $C^{k,\alpha}$ in $A_{\mu}$.
\end{theorem}

\begin{proof}
From Lemma~\ref{l:taupos2} there exists $\mu>0$
such that $\len(x_0),\len^-(x_0) > \mu$
for every $x_0\in\partial\Omega$,
hence, thanks to (\ref{f:conj}), for every point $z\in\ A_{\mu}$ the projection
$\proj^s(z)$ is a singleton.
Let $z_0\in A_{\mu}$
and prove that $\sdist$ is of class $C^2$ in a
neighborhood of $z_0$.
Let $\proj(z_0) = \{x_0\}$,
and let $Y\colon\U\to\R^n$, $\U\subset\R^{n-1}$ open,
a local parametrization of $\partial\Omega$ in a
neighborhood of $x_0$ with $Y(0) = x_0$.
Let $\Psi\in C^1(\U \times \R)$ be the map defined in (\ref{f:psi}),
and let $t_0 = \sdist(z_0)\in (-\mu,\mu)$.
{}From Theorem~\ref{t:detpsi} we have that
$\det D\Psi(0, t_0) > 0$,
hence from the inverse mapping theorem
it follows that there exists a neighborhood
$V\subset A_{\mu}$ of $z_0$ such that the maps $y = y(x)$
and $t=t(x)$ are
of class $C^1(V)$.
Since $\sdist(x) = t(x)$ for every $x\in V$,
this proves that $\sdist$ is of class $C^1(V)$.
Moreover, by  (\ref{f:diffd}),
$D\sdist(\Psi(0,t)) = \nor(x_0)/\gauge(\nor(x_0))$
for every $t\in (0,t_0)$, hence
by the continuity of $D\sdist$ in $V$, and the fact that (\ref{f:diffdb})
holds in $(0,t_0)$, we have
\begin{equation}\label{f:sulb}
D\sdist(x_0) = \nor(x_0)/\gauge(\nor(x_0))
\qquad
(x_0\in\partial\Omega)\,.
\end{equation}
Finally, from (\ref{f:diffd})
we have that
$D\dist(x) = N(y(x)) / \gauge(N(y(x)))$, $x\in V$,
where $N = \nor\circ Y\in C^1(\U)$,
hence $D\dist\in C^1(V)$, that is
$\dist\in C^2(V)$.

\noindent
The last part of the proposition
follows from the fact that,
if $\Omega$ is of class $C^{k,\alpha}$
and $\gauge\in C^{k,\alpha}(\R^n\setminus\{0\})$,
then $\Psi\in C^{k-1,\alpha}(\U\times\R)$.
\end{proof}

\begin{remark}\label{c:sig}
As a consequence of Theorem~\ref{t:regd}, we have that
$\overline{\Sigma^s}\subset\R^n\setminus\partial\Omega$.
In particular
$\overline{\Sigma}\subset\Omega$
and
$\overline{\Sigma^-}\subset\R^n\setminus\overline{\Omega}$.
\end{remark}

We conclude this section giving an explicit representation
of the Hessian matrix of $\dist$ in $x_0\in\partial\Omega$
with respect to the principal coordinate system.
From now on,
we shall denote
$D^2\dist(x_0) = D^2\sdist(x_0)$
for every $x_0\in\partial\Omega$.

\begin{lemma}\label{l:hessd}
Let $x_0\in\partial\Omega$.
Then, in the principal coordinate system
at $x_0$ we have that
\begin{gather}
\pscal{D^2\dist(x_0) e_i}{e_j} = -\frac{1}{\gauge(\nor(x_0))}\,\curv_i\,\delta_{ij},
\quad
i,j=1,\ldots,n-1.
\label{f:d2}\\
\pscal{D^2\dist(x_0) e_i}{\nor(x_0)} =
\frac{\curv_i}{\gauge(\nor(x_0))^2}
\pscal{D\gauge(\nor(x_0))}{e_i},
\
i=1,\ldots,n-1.
\label{f:d3}\\
\pscal{D^2\dist(x_0) \nor(x_0)}{\nor(x_0)} =
-\frac{1}{\gauge(\nor(x_0))^3}
\sum_{i=1}^{n-1}\curv_i\pscal{D\gauge(\nor(x_0))}{e_i}^2.
\label{f:d4}
\end{gather}
(Recall that $\nor(x_0) = e_n$ in the principal coordinate system.)
\end{lemma}

\begin{proof}
If we differentiate the identity $\dist(X(y)) = 0$ for every $y\in\U$
with respect to $y_i$, $i=1,\ldots,n-1$, we obtain
\[
\pscal{D\dist(X(y))}{\frac{\partial X}{\partial y_i}(y)} = 0,
\quad y\in\U\,.
\]
A further differentiation with respect to $y_j$,
$j=1,\ldots,n-1$, gives
\[
\pscal{D^2\dist(X(y))\frac{\partial X}{\partial y_i}(y)}%
{\frac{\partial X}{\partial y_j}(y)}
+ \pscal{D\dist(X(y))}%
{\frac{\partial^2 X}{\partial y_i\partial y_j}(y)}
= 0,
\quad y\in\U\,.
\]
Evaluating the last expression at $y=0$,
and taking into account (\ref{f:X}) and Lem\-ma~\ref{l:pis},
we obtain~(\ref{f:d2}).

\noindent
In order to prove (\ref{f:d3}) and (\ref{f:d4}),
let us differentiate the identity (\ref{f:diffdb})
with respect to $t$, and evaluate the result at $t=0$.
We obtain
\begin{equation}\label{f:eigD}
D^2\dist(x_0)\, D\gauge(\nor(x_0)) = 0,
\end{equation}
that, upon observing that
\[
D\gauge(\nor(x_0)) =
\sum_{l=1}^{n-1}\pscal{D\gauge(\nor(x_0))}{e_l}e_l
+ \gauge(\nor(x_0))\, e_n,
\]
becomes
\[
\sum_{l=1}^{n-1}
\pscal{D\gauge(\nor(x_0))}{e_l}
D^2\dist(x_0)\, e_l + \gauge(\nor(x_0))\, D^2\dist(x_0)\, e_n = 0.
\]
Then
\begin{equation}\label{f:chec}
\begin{split}
\sum_{l=1}^{n-1}&
\pscal{D\gauge(\nor(x_0))}{e_l}
\pscal{D^2\dist(x_0)\, e_l}{e_i}
\\ & +
\gauge(\nor(x_0))\, \pscal{D^2\dist(x_0)\, e_n}{e_i} = 0,
\end{split}
\end{equation}
for every $i=1,\ldots,n$. Hence
(\ref{f:d3}) follows from (\ref{f:chec}) and (\ref{f:d2}),
while (\ref{f:d4}) follows from (\ref{f:chec}) and
(\ref{f:d3}).
\end{proof}

\begin{remark}
Formula (\ref{f:d2}) can be written as
\[
\pscal{D^2\dist(x_0) e_i}{e_j} =
\frac{1}{\gauge(\nor(x_0))}\,
\pscal{D\, e_i}{e_j}\,,
\quad
i,j=1,\ldots,n-1,
\]
where $D\in\matr_{n-1}$ is the matrix
defined in (\ref{f:DR}).
\end{remark}

\section{$\gauge$-curvatures}

Thanks to Theorem~\ref{t:regd}, the function $\sdist$
is of class $C^2$ on
a neighborhood of $\partial\Omega$.
We can then define the matrix--valued function
\begin{equation}\label{f:W}
W(x_0) = -D^2\gauge(D\dist(x_0))\, D^2\dist(x_0)\,,\qquad
x_0\in\partial\Omega,
\end{equation}
where $D^2\dist(x_0) \equiv D^2\sdist(x_0)$.
For any $x_0\in\partial\Omega$ let $T_{x_0}$ denote
the tangent space to $\partial\Omega$ at $x_0$.
Notice that, if $v\in \R^n$,
then $W(x_0)\, v\in T_{x_0}$, since
from (\ref{f:sulb}) and Theorem~\ref{t:sch}(i),(iv)
\begin{equation*}
\begin{split}
\pscal{W(x_0)v}{\nor(x_0)}& =-\pscal
{D^2\dist(x_0)v}{D^2\gauge(D\dist(x_0))\, \nor(x_0)} \\
&=
-\gauge(\nor(x_0)) \pscal{D^2\dist(x_0)v}{D^2\gauge(\nor(x_0))\, \nor(x_0)} =
0.
\end{split}
\end{equation*}
Hence, we can define the map
\begin{equation}\label{f:bw}
\bw(x_0)\colon T_{x_0}\to T_{x_0},\quad
\bw(x_0)\, w = W(x_0)\, w,\quad w\in T_{x_0}\,,
\end{equation}
that can be identified with a linear application
from $\R^{n-1}$ to $\R^{n-1}$.
In our setting, the function $\bw$ plays the role of
the Weingarten map
(see Example~\ref{e:eucl} and Remark~\ref{r:wein} below).

\begin{lemma}\label{l:W1}
Let $x_0\in\partial\Omega$.
Then, in the principal coordinate system at $x_0$,
$\bw(x_0)=RD$ where $D$ and $R$ are the symmetric
matrices defined in (\ref{f:DR}),
one has
\[
\det [I-t\, W(x_0)] = \det[I_{n-1}-t\,\bw(x_0)],
\qquad t\in\R.
\]
Furthermore, both determinants are strictly positive
for $t\in (-\len^-(x_0),\len(x_0))$.
\end{lemma}

\begin{proof}
Consider the principal coordinate system at $x_0$.
Recall that, from the positive $(-1)$-homogeneity
of $D^2\gauge$ and (\ref{f:sulb}),
\[
D^2\gauge(D\dist(x_0)) = \gauge(\nor)\, D^2\gauge(\nor),
\]
where $\nor = \nor(x_0)$.
Since
$D^2\gauge(\nor)\, e_n = D^2\gauge(\nor)\, \nor = 0$
(see (\ref{f:eigenzero})),
we get
\begin{equation}\label{f:wij}
\pscal{W(x_0)\, e_j}{e_i} =
-\gauge(\nor)\sum_{l=1}^{n-1}
\pscal{D^2\gauge(\nor)\, e_l}{e_i}
\pscal{D^2\dist(x_0)\, e_j}{e_l},
\end{equation}
for every $i,j=1,\ldots,n-1$.
{}From Lemma~\ref{l:hessd} we obtain
\[
\pscal{W(x_0)\, e_j}{e_i} =
(R\, D)_{ij},
\qquad
i,j=1,\ldots,n-1,
\]
where $D$ and $R$ are the symmetric matrices
in $\matr_{n-1}$ defined in (\ref{f:DR}).
Since $T_{x_0}$ is spanned by $(e_1,\ldots,e_{n-1})$,
we have that
\begin{equation}\label{f:bwRD}
\bw(x_0) = R\, D.
\end{equation}
Furthermore
\[
\pscal{W(x_0)\, e_k}{\nor} =
-\pscal{D^2\gauge(\nor)\, \nor}{D^2\dist(x_0)\, e_k} = 0,
\qquad
k=1,\ldots,n,
\]
hence
\[
\det[I-t\, W(x_0)]
= \det[I_{n-1}-t\, \bw(x_0)]
= \det(I_{n-1} -t\,  R\, D),
\]
and the conclusion follows from
Lemmas~\ref{l:detpsi}, \ref{l:detpsi2}
and Theorem~\ref{t:detpsi}.
\end{proof}

\begin{remark}\label{r:zeig}
Let $x_0\in\partial\Omega$.
{}From (\ref{f:eigD}) we have that
\[
W(x_0)D\gauge(\nor(x_0))=
-D^2\gauge(D\dist(x_0))D^2\dist(x_0)D\gauge(\nor(x_0))=0
\]
that is
$D\gauge(\nor(x_0))$ (and hence $D\gauge(D\dist(x_0))$, by (\ref{f:sulb}))
is an eigenvector of
$W(x_0)$ with corresponding eigenvalue zero.
Since $\pscal{D\gauge(\nor(x_0))}{\nor(x_0)}=\gauge(\nor(x_0))\neq 0$,
then $D\gauge(\nor(x_0))\not\in T_{x_0}$ is not an eigenvector for
$\bw(x_0)$.
On the other hand, from Lemma~\ref{l:W1}
we deduce that a number $\curv\neq 0$ is an
eigenvalue of $W(x_0)$ if and only if
it is an eigenvalue of $\bw(x_0)$.
\end{remark}

\begin{remark}\label{r:eigv}
Although
the matrix $\bw(x_0)$ is not in general symmetric,
its eigenvalues
are real numbers
(and so its eigenvectors are real).
This property easily follows from (\ref{f:bwRD}),
and the fact that $R$ and $D$ are symmetric matrices.
Moreover, since $\bw(x)$ is a continuous matrix--valued
function on $\partial\Omega$, its eigenvalues are
continuous real functions on $\partial\Omega$.
\end{remark}

\begin{lemma}\label{l:W2}
Let $x\in\overline{\Omega}$ and let $x_0\in\proj(x)$.
If $\curv$ is an eigenvalue of $\bw(x_0)$,
then
$\curv\, \dist(x)\leq 1$.
Furthermore, if
$x\in\overline{\Omega}\setminus{\ridge}$,
then
$\curv\, \dist(x) < 1$.
\end{lemma}

\begin{proof}
If $\curv\leq 0$, then there is nothing to prove.
Assume that $\curv > 0$.
Then $\curv$ is an eigenvalue of $\bw(x_0)$ if
and only if
\[
\det[\curv\, I_{n-1}-\bw(x_0)] =
\curv^{n-1}\, \det\left[I_{n-1}-\frac{1}{\curv}\,\bw(x_0)\right] =
0\,.
\]
{}From Lemma~\ref{l:W1} this can happen
only if $1/\curv\geq\len(x_0)$,
that is, if $\curv\,\len(x_0)\leq 1$.
Since $\len(x_0)\geq \dist(x)$, this proves that
$\curv\,\dist(x)\leq 1$.
Furthermore,
if $x\in\overline{\Omega}\setminus{\ridge}$
then $\len(x_0) > \dist(x)$, so that
$\curv\,\dist(x) < 1$.
\end{proof}

The eigenvalues of $\bw(x_0)$
have a significant geometric interpretation
(see Remarks \ref{r:curvatures} and \ref{r:curvmax}).

\begin{dhef}[$\gauge$-curvatures]\label{d:curv}
Let $x_0\in\partial\Omega$.
The \textsl{principal $\gauge$-curvatures} of $\partial\Omega$ at $x_0$,
with respect to the Minkowski norm $\dist$,
are the eigenvalues
$\curvg_1(x_0)\leq\cdots\leq\curvg_{n-1}(x_0)$
of $\bw(x_0)$.
The corresponding eigenvectors are the
\textsl{principal $\gauge$-directions} of $\partial\Omega$
at $x_0$.
\end{dhef}

\begin{example}[Euclidean distance]\label{e:eucl}
Let $K = \overline{B}_1(0)$,
so that $\gauge(\xi) = |\xi|$.
In this case, $\dist$ coincides with
the Euclidean distance $\dist^E$ from
$\partial\Omega$.
Let $x_0\in\partial\Omega$
and consider the principal coordinate system at $x_0$.
Since
\[
D\gauge(\xi) = \frac{\xi}{|\xi|}\,,\quad
D^2\gauge(\xi) = \frac{1}{|\xi|}\, I -
\frac{1}{|\xi|^3}\,\xi\otimes\xi,
\quad\xi\neq 0,
\]
{}from
Lemma~\ref{l:hessd} we recover
$W(x_0) = \diag(\curv_1,\ldots,\curv_{n-1}, 0)=-D^2\dist^E(x_0)$ and
$\bw(x_0) = \diag(\curv_1,\ldots,\curv_{n-1})$.
Hence the principal $\gauge$-curvatures and the
principal $\gauge$-directions correspond
respectively to the principal curvatures and
the principal directions of $\partial\Omega$ at $x_0$.
\end{example}

\begin{remark}[Normal curvatures]\label{r:curvatures}
Definition~\ref{d:curv} is motivated by the following construction.
Let $x_0\in\partial\Omega$ and consider the
principal coordinate system at $x_0$.
Let $X\colon\U\to\R^n$ be a standard parametrization of $\partial\Omega$
in a neighborhood of $x_0$ satisfying (\ref{f:X}).
Let $v\in\R^{n-1}$, $v\neq 0$, and consider the curve on $\partial\Omega$
\[
x(t) = X(tv),\quad
t\in I,
\]
where $I = (-t_0,t_0)$ and $t_0>0$ is chosen such
that $t\,v\in\U$ for every $t\in I$.
Let us denote by $V$ the plane generated by $v$ and $\nor=\nor(x_0) = e_n$.
For $r>0$, let $z_r = x_0 + r\, D\gauge(\nor)$
and let $K_r = z_r - r\, K^0$, so that
$x_0\in K_r\cap\partial\Omega$.
We would like to determine $r>0$ in such a way that
the section $K_r\cap V$ has a contact of order two
with the curve $x(t)$ at $x_0$.
Such a value of $r=r(v)$ will be the radius of curvature
of the curve $x(t)$ at $x_0$
with respect to the Minkowski distance, hence
$\curvg(v) = 1/r(v)$
will be the curvature of $x(t)$ at $x_0$,
that is, the normal curvature of $\partial\Omega$ at $x_0$
in the direction of $v$ with respect to the
Minkowski distance.

\noindent
For $t\in I$,
the distance (with respect to the Minkowski norm)
from $x(t)$ to the section $K_r\cap V$ can be
estimated by $|\pgauge(z_r-x(t))-r|$.
The condition of a contact of second order becomes
\begin{equation}\label{f:derh}
\begin{split}
0 & = \left.\frac{d^2}{dt^2}\,
\pgauge(z_r-x(t))\right|_{t=0}
\\ & =
\pscal{D^2\pgauge(z_r-x_0)\,\dot{x}(0)}{\dot{x}(0)}+
\pscal{D\pgauge(z_r-x_0)}{\ddot{x}(0)}\,.
\end{split}
\end{equation}
Now, $z_r-x_0 = r\, D\gauge(\nor)$,
so that from homogeneity and Lemma~\ref{l:polare}
we get
\[
D^2\pgauge(z_r-x_0) = \frac{1}{r}\, D^2\pgauge(D\gauge(\nor)),
\quad
D\pgauge(z_r-x_0) = D\pgauge(D\gauge(\nor)) = \frac{\nor}{\gauge(\nor)}.
\]
On the other hand,
using (\ref{f:fiz}) and (\ref{f:X}) we get
\[
\dot{x}(0) = (v,0),\qquad
\ddot{x}(0) = \left(0, \sum_{i=1}^{n-1}\curv_i\, v_i^2\right),
\]
hence (\ref{f:derh}) becomes
\begin{equation}\label{f:derh2}
\frac{1}{r}\,\sum_{i,j=1}^{n-1}
\pscal{D^2\pgauge(D\gauge(\nor))\, e_j}{e_i}\, v_i\, v_j
+ \frac{1}{\gauge(\nor)}\,\sum_{i=1}^{n-1}\curv_i\, v_i^2
= 0\,.
\end{equation}
Recalling the definitions of the matrices $H$ and $D$
in (\ref{f:H}) and (\ref{f:DR}) respectively,
we conclude that
\[
\frac{1}{r}\, \pscal{H\, v}{v} -
\frac{1}{\gauge(\nor)}\,\pscal{D\, v}{v} = 0\,.
\]
Let us define $\curvg = \curvg(v) = 1/r$.
The identity above can we rewritten as
\begin{equation}\label{f:curvn}
\curvg(v)\, \pscal{\gauge(\nor)\,H\, v}{v} -
\pscal{D\, v}{v} = 0\,.
\end{equation}
By construction, the quantity $\curvg(v)$
is the normal curvature of $\partial\Omega$ at $x_0$
along the direction $v$, with respect to the Minkowski norm.
In analogy of what happens in the Riemannian case,
the principal $\gauge$-curvatures
can be defined as the invariants
of the pair of quadratic forms $D$
and $\gauge(\nor)\, H$,
that is, the numbers $\lambda\in\R$
that satisfy the characteristic equation
\[
\det(D-\lambda\,\gauge(\nor)\,H) = 0\,.
\]
On the other hand,
from (\ref{f:RH}) we have that $(\gauge(\nor)\, H)^{-1} = R$,
where $R$ is the matrix defined in (\ref{f:DR}),
hence the invariants are the solutions of
$\det(R\, D-\lambda\,I_{n-1}) = 0$.
Recalling (\ref{f:bwRD}), these are the solutions to
\[
\det(\bw(x_0)-\lambda\,I_{n-1}) = 0\,,
\]
that is, the eigenvalues of $\bw(x_0)$.
Moreover, if
$\curvg_1(x_0)\leq\cdots\leq\curvg_{n-1}(x_0)$
are the eigenvalues of $\bw(x_0)$ and
$v_1,\ldots,v_{n-1}$ are the corresponding eigenvectors,
then it is straightforward to check that
$\curvg_i(x_0) = \curvg(v_i)$ for every
$i=1,\ldots,n-1$.
\end{remark}

\begin{remark}\label{r:curvmax}
Let $x_0\in\partial\Omega$ and let
$\curvg_1(x_0)\leq\cdots\leq\curvg_{n-1}(x_0)$
be the principal $\gauge$-curvatures of $\partial\Omega$
at $x_0$, that is, the eigenvalues of $\bw(x_0)$.
We claim that
\begin{equation}\label{f:norcurv}
\curvg_1(x_0)\leq\curvg(v)\leq\curvg_{n-1}(x_0)
\qquad
\forall v\in U_{x_0} :=
\{w\in T_{x_0}\Omega;\
\|w\| = 1\}\,,
\end{equation}
where $\curvg(v)$ is the normal curvature of $\partial\Omega$
at $x_0$ along the direction $v$,
defined in (\ref{f:curvn}),
and $U_{x_0}$ is the unit tangent bundle of $\partial\Omega$
at $x_0$.
More precisely, we claim that
\begin{equation}\label{f:minmax}
\curvg_1(x_0) =
\min_{v\in U_{x_0}}\curvg(v),\qquad
\curvg_{n-1}(x_0) =
\max_{v\in U_{x_0}}\curvg(v)\,.
\end{equation}
Namely, in the principal coordinate system at $x_0$,
$\curvg\colon\R^{n-1}\setminus\{0\}\to\R$
is a continuous function on $S^{n-2}$,
hence it admits maximum and minimum on $S^{n-2}$.
Let $\bar{v}\in S^{n-2}$ be a maximum point.
Since $\curvg$ is positively $0$-homogeneous (see (\ref{f:curvn})),
then $\bar{v}$ is also a maximum point of $\curvg$
in the open set $\R^n\setminus\{0\}$,
so that $D\curvg(\bar{v}) = 0$.
Differentiating (\ref{f:curvn}) and plugging this
identity we obtain
\[
\curvg(\bar{v})\, \gauge(\nor)\,H\, \bar{v} - D\, \bar{v} = 0\,,
\]
that is,
$\det(\curvg(\bar{v})\, \gauge(\nor)\,H - D) = 0$.
Thus $\curvg(\bar{v})$ is an invariant of the
pair of quadratic forms $D$ and $\gauge(\nor)\, H$,
and hence an eigenvalue of $\bw(x_0)$.
Finally, being $\curvg(\bar{v})$ the maximum of $\curvg$
on $S^{n-2}$, we conclude that
$\curvg(\bar{v}) = \curvg_{n-1}(x_0)$.
Reasoning as above, if $\bar{w}$ is a minimum point
of $\curvg$ on $S^{n-2}$, we deduce that
$\curvg(\bar{w}) = \curvg_{1}(x_0)$.
\end{remark}

\begin{remark}[Anisotropic Weingarten map]\label{r:wein}
We can also give another interpretation of
$\gauge$-curvatures (see \cite[\S2.3]{An}).
We recall that, by definition, the principal curvatures
of a smooth manifold $M$ at a given point $x\in M$
are the eigenvalues of the Weingarten map
$L_x = -d\nor_x\colon T_x M\to T_x M$,
where $\nor(x)$ is the normal to $M$
at $x$.
In out setting, the inward ``normal'' to
$\partial\Omega$ at $x$ is given by the
vector $p(x):=D\gauge(\nor(x))$
(see Proposition~\ref{p:pnc} and Remark~\ref{r:pnc}).
Hence the Weingarten map, in the Minkowskian setting,
should be $L^M_x = -dp_x\colon T_x\Omega\to T_x\Omega$.
It is readily seen that $L^M_x$ coincides with the
map $\bw(x)$ defined in (\ref{f:bw}).
Hence, the principal curvatures at $x$ with respect to
the Minkowski distance are the eigenvalues
of the anisotropic Weingarten map $\bw(x)$.
\end{remark}

\begin{remark}[Curvatures in Finsler geometry]
If $\pgauge$ is of class $C^{\infty}$,
then the principal $\gauge$-curvatures
coincide with the notion of curvature
in the setting of Finsler geometry
as defined in \cite{Sh}.
\end{remark}

\section{Regularity of $\overline{\Sigma}$}

In this section we investigate the structure of the set
$\overline{\Sigma}$. We shall prove that $\overline{\Sigma}=
\ridge$, and that $\dist$ is a function of class $C^2$ in
the whole $\overline{\Omega}\setminus\overline{\Sigma}$.
Finally we propose a regularity result for the set
$\overline{\Sigma}$.
The main tool for the proof of these results will be a complete
description of the set $\overline{\Sigma}\setminus\Sigma$.

\begin{dhef}[Focal points]
Let us define the map
\begin{equation}\label{f:Phi}
\Phi(x,t) = x+t\, D\gauge(\nor(x)),
\qquad x\in\partial\Omega,\ t\in\R.
\end{equation}
Fixed $x_0\in\partial\Omega$,
we say that $z=x_0+t_0\,D\gauge(\nor(x_0))$
is a focal point of $x_0$
if the map $\Phi$ is singular at $(x_0, t_0)$, that is
$\det D\Phi(x_0,t_0)=0$.
\end{dhef}

In the setting of Riemannian geometry,
the point $z$ in the above definition
would be called a focal point of
$\partial\Omega$ along the geodesic
$t\mapsto \Phi(x_0, t)$, $t\geq 0$,
normal to $\partial\Omega$ at $x_0$
(see e.g.~\cite{Sa}).

\begin{remark}\label{r:llessf}
Let $Y\colon\U\to\R^n$,
$\U\subset\R^{n-1}$ open,
be a local parametrization
of $\partial\Omega$ in
a neighborhood of $x_0$, such that $Y(0) = x_0$.
The fact that $z=x_0+t_0 D\gauge(\nor(x_0))$ is a
focal point amounts to say that the differential of the map
$\Psi(y,t) := Y(y)+t\, D\gauge(\nor(Y(y)))$
at $(0, t_0)$ has not maximal rank, that is
$\det D\Psi(0, t_0) = 0$.
Recalling (\ref{f:calcdetpsi}), (\ref{f:calcdetpsi2}) and
(\ref{f:bwRD}), we have
\[
\det D\Psi(0, t_0) = 0
\qquad\Longleftrightarrow\qquad
\det(I_{n-1}-t_0 \, \bw(x_0)) = 0,
\]
hence the focal points of
$x_0\in\partial\Omega$ are
$\Phi(x_0, 1/\curvg_i(x_0))$,
if $\curvg_i(x_0)\neq 0$,
$i=1,\ldots,n-1$. Moreover, by Theorem~\ref{t:detpsi},
$\len(x_0)\leq 1/\curvg_{n-1}(x_0)$.
\end{remark}

\begin{dhef}[Optimal focal points]\label{d:focal}
Let $x_0\in\partial\Omega$.
If $\curvg_{n-1}(x_0) > 0$,
we call $z_0 = x_0 + 1/\curvg_{n-1}(x_0)\, D\gauge(\nor(x_0))$
the first focal point of $x_0$.
If, in addition, $z_0\in\overline{\Sigma}$,
we call $z_0$ the optimal focal point of $x_0$.
We denote by $\Gamma$ the set of optimal focal points
of the boundary points of $\Omega$.
\end{dhef}

In the following proposition we prove that
$\overline{\Sigma}= \Sigma\cup\Gamma$.

\begin{prop}\label{p:noconj}
If $x\in\overline{\Sigma}\setminus\Sigma$,
then
$\curvg_{n-1}(x_0)\, \dist(x) = 1$,
where $x_0$ is the unique projection of $x$
on $\partial\Omega$.
\end{prop}

\begin{proof}
Let $x\in\overline{\Sigma}\setminus\Sigma$, and
let $(x_k)\subset\Sigma$ be a sequence converging
to $x$.
For every $k\in\N$, let $y_k, z_k\in\partial\Omega$
be two distinct points in $\proj(x_k)$,
so that
\begin{equation}\label{f:xautov}
x_k = y_k+\dist(x_k)\, D\gauge(D\dist(y_k))
= z_k+\dist(x_k)\, D\gauge(D\dist(z_k))\,,
\end{equation}
and let
\[
\xi_k = \frac{y_k-z_k}{|y_k-z_k|}\,.
\]
Since $\xi_k\in S^{n-1}$ for every $k$,
we can extract a subsequence (which we do not relabel)
converging to a point $\xi\in S^{n-1}$.
From Lemma \ref{l:pusc} we have that
$(y_k)_k$ and $(z_k)_k$ converge to $x_0$.
From (\ref{f:xautov}) we get
\[
\frac{y_k-z_k}{|y_k-z_k|}
+\dist(x_k)\, \frac{D\gauge(D\dist(y_k))-D\gauge(D\dist(z_k))}{|y_k-z_k|}
= 0\quad
\forall k\in\N.
\]
Since $\dist$ is of class $C^2$ in a neighborhood of $\partial\Omega$
and $\gauge$ is of class $C^2$ in $\R^n\setminus \{0\}$,
passing to the limit we get
$\xi + \dist(x)\, D^2\gauge(D\dist(x_0))\, D^2\dist(x_0)\, \xi = 0$,
that is
$\dist(x)\, W(x_0)\, \xi = \xi$.
Hence $1/\dist(x)$ is a nonzero eigenvalue of $W(x_0)$.
From Remark~\ref{r:zeig}, we conclude that
$1/\dist(x)$ is an eigenvalue of $\bw(x_0)$, that is,
$1/\dist(x) = \curvg_j(x_0)$ for some $j\in\{1,\ldots,n-1\}$.

Finally, let us prove that
$\curvg_j(x_0) = \curvg_{n-1}(x_0)$.
Assume, by contradiction, that
$\curvg_j(x_0) < \curvg_{n-1}(x_0)$.
Then
$1/\curvg_{n-1}(x_0) < 1/\curvg_j(x_0)=\dist(x_0)\leq \len(x_0)$,
a contradiction (see Remark \ref{r:llessf}).
\end{proof}

\begin{prop}\label{p:cutray}
Let $x\in\Omega$ and $x_0\in\proj(x)$.
Then $\overline{\Sigma}\cap (x_0,x) = \emptyset$.
\end{prop}

\begin{proof}
By Proposition~\ref{p:basic}(iii) we already know that
$\Sigma\cap (x_0, x) = \emptyset$.
Moreover, from Lemma~\ref{l:W2}
we have that $\curvg_{n-1}(x_0)\, \dist(x)\leq 1$,
hence
the strict inequality
$\curvg_{n-1}(x_0)\, \dist(z) < 1$
holds for every $z\in (x_0,x)$.
This implies, by Proposition \ref{p:noconj}, that
$(\overline{\Sigma}\setminus\Sigma)\cap (x_0, x) = \emptyset$,
concluding the proof.
\end{proof}

\begin{cor}\label{c:bang}
For every $x_0\in\partial\Omega$,
\begin{equation}\label{f:lenmin}
\len(x_0) =
\min\{t\geq 0;\
x_0+t\, D\gauge(\nor(x_0))\in\overline{\Sigma}\}\,.
\end{equation}
\end{cor}

\begin{proof}
Let us denote by $\lambda(x_0)$ the
right--hand side of (\ref{f:lenmin}).
Notice that the point
$z_0=x_0+\len(x_0) D\gauge(\nor(x_0))$
belongs to the ridge set, then, by
Remark~\ref{r:bors}, $z_0\in \overline{\Sigma}$.
Hence we have $\lambda(x_0) \leq \len(x_0)$.
Assume by contradiction that
the strict inequality holds.
Then, if we set
$z(t) := x_0+t\, D\gauge(\nor(x_0))$,
for $\overline{t}\in (\lambda(x_0), \len(x_0))$
we would have $x_0\in\proj(z(\overline{t}))$ and
$z(\lambda(x_0))\in (x_0, z(\overline{t}))\cap\overline{\Sigma}$,
in contradiction with Proposition~\ref{p:cutray}.
\end{proof}

\begin{theorem}\label{t:tau}
The function $\len$, defined in (\ref{f:len}),
and the function $\tau\colon\overline{\Omega} \to \R$
defined by
\begin{equation}\label{f:tau}
\tau(x) =
\begin{cases}
\min\{t\geq 0;\
x + t D\gauge(D\dist(x))\in\overline{\Sigma}\},
&\textrm{if $x\in\overline{\Omega}\setminus\overline{\Sigma}$},\\
0,
&\textrm{if $x\in\overline{\Sigma}$}\,,
\end{cases}
\end{equation}
(and called  normal distance to cut locus of $x$)
are continuous in $\overline{\Omega}$.
\end{theorem}

\begin{proof}
It is easy to check that,
if $x\in\overline{\Omega}$
and $x_0\in\proj(x)$,
then
$\tau(x) = \tau(x_0) - \dist(x)$.
Since, by Corollary~\ref{c:bang}, $\tau(x_0) = \len(x_0) = \len(x)$,
we get $\tau(x) = \len(x)-\dist(x)$
for every
$x\in\overline{\Omega}$.
From Proposition~\ref{p:lusc} the
function $\len$ is upper semicontinuous in $\Omega$,
hence the upper semicontinuity of $\tau$ follows from
the continuity of $\dist$.
On the other hand,
the lower semicontinuity of $\tau$
at a point $x\in\overline{\Sigma}$ is trivial.
Hence it remains to prove that $\tau$ is lower semicontinuous
in $x\in\overline{\Omega}\setminus\overline{\Sigma}$.
Let $(x_k)_k$ be a sequence in $\overline{\Omega}$ converging to $x$.
We can assume that $x_k\in \overline{\Omega}\setminus\overline{\Sigma}$,
hence, by definition of $\tau$,
$x_k + \tau(x_k)\, D\gauge(D\dist(x_k))\in\overline{\Sigma}$.
If $\tau_0 := \liminf_k \tau(x_k)$,
passing to the limit we get
$x+\tau_0\, D\gauge(D\dist(x))\in\overline{\Sigma}$,
that is, $\tau(x)\leq \tau_0$.
\end{proof}

\begin{cor}\label{c:closed}
The ridge set $\ridge$ is closed and
$\ridge = \overline{\Sigma}$.
\end{cor}

\begin{proof}
We recall that
$\ridge = \{x\in\Omega;\ \len(x) = \dist(x)\}$.
Since, by Lemma~\ref{l:taupos2}, $\overline{\ridge}\subset\Omega$,
the fact that $\ridge$ is closed follows
from the continuity of $\len$ and $\dist$.
The equality $\ridge = \overline{\Sigma}$
now follows from the inclusions $\Sigma\subseteq \ridge\subseteq
\overline{\Sigma}$.
\end{proof}

\begin{cor}\label{c:mn}
The set $\overline{\Sigma}$ has
vanishing Lebesgue measure.
\end{cor}

\begin{proof}
It follows from
Corollaries~\ref{c:sigma}
and \ref{c:closed}.
\end{proof}

We can now prove a regularity result for $\dist$,
which extends Theorem~\ref{t:regd}
to $\overline{\Omega}\setminus\overline{\Sigma}$.

\begin{theorem}[Regularity of $\dist$ in
$\overline{\Omega}\setminus\overline{\Sigma}$]\label{t:regd2}
The function $\dist$ is of class $C^2$ in
$\overline{\Omega}\setminus\overline{\Sigma}$.
If
$\Omega$ is of class $C^{k,\alpha}$
and $\gauge\in C^{k,\alpha}(\R^n\setminus\{0\})$,
for some $k\geq 2$ and $\alpha\in [0,1]$,
then $\dist$ is of class $C^{k,\alpha}$ in
$\overline{\Omega}\setminus\overline{\Sigma}$.
\end{theorem}

\begin{proof}
It is enough to
prove the result in $\Omega\setminus\overline{\Sigma}$, since
the regularity of $\dist$ near $\partial \Omega$ was already
proved in Theorem~\ref{t:regd}.

Fixed $z_0\in\Omega\setminus\overline{\Sigma}$,
we have to show that $\dist$ is of class $C^2$ in a
neighborhood of $z_0$.
Let $\proj(z_0)=\{x_0\}$,
and let $Y\colon\U\to\R^n$, $\U\subset\R^{n-1}$ open,
be a local parametrization of $\partial\Omega$ in a
neighborhood of $x_0$.
Let $\Psi$ be the map defined in (\ref{f:psi}),
and let $t_0 = \dist(z_0)$. By Corollary~\ref{c:closed},
$z_0$ does not belong to $\ridge$, that is
$\dist(z_0) < \len(z_0)=\len(x_0)$.
{}From Theorem~\ref{t:detpsi} we have that
$\det D\Psi(x_0, t_0) > 0$,
hence from the inverse mapping theorem
it follows that there exists a neighborhood
$V$ of $z_0$ such that the map $y = y(x)$ is
of class $C^1(V)$.
On the other hand, from Lemma~\ref{l:pis},
$D\dist(x) = N(y(x)) / \gauge(N(y(x)))$, $x\in V$,
where $N = \nor\circ Y$, hence
we conclude that $D\dist\in C^1(V)$, that is
$\dist\in C^2(V)$.

\noindent
The last part of the proposition
follows from the fact that,
if $\Omega$ is of class $C^{k,\alpha}$
and $\gauge\in C^{k,\alpha}(\R^n\setminus\{0\})$,
then $\Psi\in C^{k-1,\alpha}(\U\times\R)$.
\end{proof}

In literature the set
\[
\cutl(\Omega)=\{\Phi(x_0,\len(x_0)) = x_0 + \len(x_0)\, D\gauge(\nor(x_0)),\
x_0\in\partial\Omega\}
\]
is called
the \textsl{cut locus} of $\Omega$.
From (\ref{f:lenmin}), Corollary~\ref{c:closed},
and Proposition~\ref{p:noconj}
it follows that
\[
\cutl(\Omega) = \ridge = \overline{\Sigma} = \Sigma\cup\Gamma.
\]
The regularity of the singular set $\Sigma$
and of the cut locus $\overline{\Sigma}$ has been
extensively studied.
For what concerns $\Sigma$,
Alberti \cite{Alb} has proved that
it is $C^2$-rectifiable.
We recall that
a subset of $\R^n$ is $C^k$-rectifiable, $k\in\N$,
if it can be covered by a countable family of
embedded $C^k$ manifolds of dimension $n-1$,
with the exception of a set of vanishing $\haus$ measure.
On the other hand, in general
the cut locus $\overline{\Sigma}$ is not rectifiable.
Namely,
Mantegazza and Mennucci \cite{MM}
have exhibited a set $\Omega\subset\R^2$ of class
$C^{1,1}$ such that the singular set $\Sigma$ corresponding
to the Euclidean distance $\dist^E$
has closure $\overline{\Sigma}$
with positive Lebesgue measure.
In the same paper the authors have proven
that, if $\Omega$ is an open subset
of class $C^r$, with $r\geq 3$,
of a smooth, connected
and complete Riemannian
manifold without boundary,
then $\cutl(\Omega)$ is $(r-2)$-rectifiable.
In the same setting,
Itoh and Tanaka \cite{IT} have proven that,
if $\Omega$ is of class $C^{\infty}$,
then the function $\len$ is Lipschitz continuous on
$\partial\Omega$.
Recently, Li and Nirenberg \cite{LN}
have refined this result, extending it
to the case of Finsler manifolds.
Adapting their result to our setting,
they have proved that,
if $\Omega$ is of class $C^{2,1}$
and $\gauge$ is of class $C^{\infty}$,
then $\len$ is Lipschitz continuous on $\partial\Omega$
(see \cite{LN}, Theorem~1.5).
As a straightforward consequence of the fact that $\len$ is
Lipschitz continuous on $\partial\Omega$, one has that
the cut locus $\overline{\Sigma}$ has finite
$\haus$ measure and it is $C^1$-rectifiable.
Namely, if $L$ denotes the Lipschitz constant
of $\len$ in $\partial\Omega$, we have
$\haus(\overline{\Sigma})\leq L\, \haus(\partial\Omega)
< +\infty$.
Other results in this direction are proved in
\cite{Me},
again in the setting of Finsler manifolds.
We mention also \cite{CaMS, CaSi} for further regularity results for
Hamilton-Jacobi equations, and \cite{Pig} for rectifiability
results of $\overline{\Sigma}$ from an optimal control theory
viewpoint.

\smallskip
In the following theorem we shall prove a
$C^{0,\alpha}$-rectifiability result of
the cut locus in the Minkowskian setting.
For what concerns the regularity of the cut locus,
our result is finer
than the ones given in \cite{LN, Me},
that, on the other hand,
deal with a general Finsler manifold.

\begin{theorem}[Regularity of $\overline{\Sigma}$]\label{t:regul}
Assume that,
for some $\alpha\in [0,1]$,
$\Omega$ is of class
$C^{2,\alpha}$,
and
$\gauge,\pgauge\in C^{2,\alpha}(\R^n\setminus\{0\})$.
Then $\overline{\Sigma}$ can be covered by countably many
graphs of functions of class $C^{0,\alpha}(\R^{n-1})$,
with the exception of a set of vanishing
$\haus$ measure.
As a consequence,
$\dimh (\overline{\Sigma})\leq n-\alpha$,
where $\dimh (\overline{\Sigma})$ denotes
the Hausdorff dimension of
$\overline{\Sigma}$.
\end{theorem}

\begin{proof}
Since $\dist$ is semiconcave in $\Omega$,
by Theorem~1 in \cite{Alb}
we have that $\Sigma$ is $C^2$-rectifiable.
Hence it is enough to prove that
the set of optimal focal points $\Gamma$
can be covered by countably many
graphs of functions of class $C^{0,\alpha}$,
defined on open subsets of $\R^{n-1}$
($C^{0,\alpha}$-rectifiable for short).
Since $\Gamma\subseteq \ridge$ is contained in
the set
\[
\Gamma_0 = \left\{
x_0+\frac{1}{\curvg_{n-1}(x_0)}\, D\gauge(\nor(x_0));\
x_0\in\partial\Omega\ \textrm{such that\ }\
\curvg_{n-1}(x_0)\geq 1/L
\right\}\,,
\]
of all first focal points
with distance from $\partial\Omega$ not
exceeding the quantity $L := \max\{\len(y);\ y\in\partial\Omega\}$,
then it is enough to prove
that $\Gamma_0$ is $C^{0,\alpha}$-rectifiable.
Let us denote by $S_0$ the set
\[
S_0 := \{x_0\in\partial\Omega;\
\curvg_{n-1}(x_0)\geq 1/L\}.
\]
Let $\U_1,\ldots,\U_N\subset\R^{n-1}$ be open sets,
and $Y_k\colon\U_k\to\R^n$, $k=1,\ldots, N$,
maps of class $C^{2,\alpha}$
that parameterize $\partial\Omega$.
For every $k=1,\ldots, N$, let
$U_k:= \{y\in\U_k;\ Y_k(y)\in S_0 \}$.
Then $S_0 = \bigcup_{k=1}^{N} Y_k(U_k)$,
and
\begin{equation}\label{f:gamma}
\Gamma_0 = \bigcup_{k=1}^{N}
\left\{
\Psi_k\left(y, \frac{1}{\curvg_{n-1}(Y_k(y))}\right);\
y\in U_k
\right\}\,,
\end{equation}
where, for every $k=1,\ldots, N$,
\[
\Psi_k(y,t) := Y_k(y)+ t\, D\gauge(\nor(Y_k(y)))\qquad
(y,t)\in\U_k\times\R\,.
\]
Let us fix $k=1,\ldots,N$.
We claim that the map
$K :=\curvg_{n-1}\circ Y_k$ is of class $C^{0,\alpha}$ in $\U_k$.
For every $y\in\U_k$ let us denote by
$D(y)$ and $H(y)$ the matrices defined respectively in
(\ref{f:DR}) and (\ref{f:H})
in the point $x_0 = Y_k(y)$.
We recall that $D(y)$ and $H(y)$ are the restrictions
respectively
of $D^2\dist(Y_k(y))$ and $D^2\pgauge(D\gauge(\nor(Y_k(y))))$
to the tangent space $T_{Y_k(y)}\Omega$.
From Theorem~\ref{t:regd2},
these functions are both of class $C^{0,\alpha}(\U_k)$.
Thus, the matrix
\[
H_0(y) := \gauge(\nor(Y_k(y)))\, H(y),
\qquad y\in\U_k,
\]
is of class $C^{0,\alpha}(\U_k)$.
Furthermore,
from Theorem~\ref{t:sch}(iv), the minimum eigenvalue
of $H(y)$ is bounded from below by the
minimum for $\xi\in\partial K$ of the radii of curvature
of $\partial K$ at $\xi$
(which is strictly positive by the assumption
$K\in C^2_+$).
From (\ref{f:brho}) we conclude that there exists
a constant $c>0$ such that
\begin{equation}\label{f:hz}
\pscal{H_0(y)\, v}{v}\geq c\|v\|^2,
\qquad
\forall y\in\U_k,\
v\in\R^{n-1}\,.
\end{equation}
From (\ref{f:minmax}),
$K(y)$ is characterized by
\begin{equation}\label{f:Kmax}
K(y) = \max_{v\in S^{n-2}}
\frac{\pscal{D(y)\, v}{v}}%
{\pscal{H_0(y)\, v}{v}}\,,
\qquad y\in\U_k.
\end{equation}
Now, let $y,y'\in\U_k$, and denote by $v_0, v_0'\in\R^{n-1}$
two vectors that realize the maximum
in (\ref{f:Kmax}) for $K(y)$ and $K(y')$ respectively.
Since $D, H_0\in C^{2,\alpha}(\U_k)$,
from (\ref{f:hz}) we have that
\[
\begin{split}
K(y)-K(y')
& \leq
\frac{\pscal{D(y)\, v_0}{v_0}}{\pscal{H_0(y)\, v_0}{v_0}} -
\frac{\pscal{D(y')\, v_0}{v_0}}%
{\pscal{H_0(y')\, v_0}{v_0}}
\\ & \leq
C\, \frac{\|H_0(y')\| + \|D(y')\| }%
{\pscal{H_0(y)\, v_0}{v_0}\,
\pscal{H_0(y')\, v_0}{v_0}}
\, |y-y'|^{\alpha}
\\ & \leq
\frac{C\, S}{c^2}\,  |y-y'|^{\alpha}
\,,
\end{split}
\]
where $S:=\sup\{\|H_0(z)\| + \|D(z)\|;\
z\in\U_k\}$.
Exchanging the role of $y$ and $y'$,
we obtain
$K(y')-K(y)\leq (C\, S/c^2)\, |y-y'|^{\alpha}$,
hence $K\in C^{0,\alpha}(\U_k)$.

\noindent
Since $\Psi_k\in C^1(\U_k\times\R)$,
with bounded derivatives in $\U_k\times [0,L]$
(see Lem\-ma~\ref{l:detpsi}),
it is Lipschitz continuous in $\U_k\times [0,L]$,
say of rank $L_k$,
hence
\[
\mathcal{H}^{n-\alpha} \{\Psi_k(y, 1/K(y));\ y\in U_k\}
\leq
L_k^{n-\alpha}\,
\mathcal{H}^{n-\alpha} \{(y, 1/K(y));\ y\in U_k\}\,.
\]
On the other hand,
being $K\colon\U_k\subset\R^{n-1}\to\R$
a positive function of class $C^{0,\alpha}$
on $\U_k$, with $K(y)\geq 1/L$ for every $y\in U_k$,
the function $1/K$ is of class $C^{0,\alpha}$
on $U_k$.
Thus the Hausdorff $(n-\alpha)$-dimensional
measure of the graph of $1/K$ on $U_k$ is finite.
Since this fact holds for every
$k=1,\ldots,N$,
from the definition (\ref{f:gamma}) of $\Gamma_0$
we conclude that
$\mathcal{H}^{n-\alpha}(\Gamma_0)$ is finite.
\end{proof}

\section{An application to PDEs}\label{s:PDE}

This section is devoted to an application
of the previous results
to the analysis of PDEs arising from
optimal transportation theory and shape optimization
(see \cite{BoBu, CCCG}).
The main tool is the following change of variables
theorem,
based on Lemma~\ref{l:detpsi},
Theorem~\ref{t:detpsi}, and the fact that
$\overline{\Sigma}$ has vanishing Lebesgue
measure.

\begin{theorem}[Change of variables]\label{t:chvar}
For every $h\in L^1(\Omega)$
\begin{equation}\label{f:chvar}
\begin{split}
& \int_{\Omega} h(x)\, dx
\\ & =
\int_{\partial\Omega} \gauge(\nor(x))\, \left[
\int_{0}^{\len(x)}
h(\Phi(x,t))\,
\det(I_{n-1}-t\, \bw(x))\, dt\right]\,
d\haus(x)
\\ & =
\int_{\partial\Omega} \gauge(\nor(x))\, \left[
\int_{0}^{\len(x)}
h(\Phi(x,t))\,
\prod_{i=1}^{n-1} (1-t\, \curvg_i(x))
\, dt\right]\,d\haus(x)
\end{split}
\end{equation}
where $\Phi\colon\partial\Omega\times\R\to\R^n$ is
the map defined by
$\Phi(x,t) = x + t\, D\gauge(\nor(x))$,
$(x,t)\in\partial\Omega\times\R$.
\end{theorem}

\begin{proof}
Let $Y_k\colon\U_k\to\R^n$, $\U_k\subset\R^{n-1}$ open,
$k=1,\ldots,N$,
be local parameterizations of $\partial\Omega$ of class $C^2$,
such that $\bigcup_{k=1}^{N} Y_k(\U_k) = \partial\Omega$.
For every $k=1,\ldots,N$,
let $\Psi_k\colon \U_k\times\R\to\R^n$ be the map
\begin{equation}\label{f:psik}
\Psi_k(y,t) = Y_k(y) + t D\gauge(\nor(Y_k(y))),
\quad (y,t)\in \U_k\times\R\,.
\end{equation}
{}From Lemma~\ref{l:detpsi}
and Theorem~\ref{t:detpsi} we know that,
for every $k=1,\ldots,N$,
$\Psi_k\in C^1(\U_k\times\R)$, and
\[
\begin{split}
\det D\Psi_k(y,t) & = \gauge(\nor(Y_k(y))\,
\sqrt{g^k(y)}\,
\det[I_{n-1}-t\, \bw(Y_k(y))]
> 0
\end{split}
\]
for every
$y\in\U_k$ and every $t\in [0, \len(Y_k(y)))$,
where
\[
g^k(y) = \det (g^k_{ij}(y)),
\quad
g^k_{ij}(y) = \pscal{\frac{\partial Y_k}{\partial y_i}(y)}%
{\frac{\partial Y_k}{\partial y_j}(y)}
\,,\quad
i,j=1,\ldots,n-1.
\]
Let $p_1,\ldots,p_N\in\partial\Omega\to\R$
be a partition of unity of $\partial\Omega$
subordinate to $Y_1,\ldots,Y_N$, that
is, for every $k=1,\ldots,N$,
$p_k(x)\geq 0$ for every $x\in\partial\Omega$,
$p_k$ has compact support contained in $Y_k(\U_k)$,
$p_k \circ Y_k\in C^{2}(\U_k)$,
and $\sum_{k=1}^N p_k(x) = 1$ for every $x\in\partial\Omega$.
Moreover,
for every $k=1,\ldots,N$
let us define the function
\[
q_k(x) =
\begin{cases}
p_k(x_0),
& \textrm{if $x_0\in Y_k(\U_k)$ and $x = \Phi(x_0,t)$
for some $t\in [0,\len(x_0))$},\\
0,
&\textrm{otherwise}\,,
\end{cases}
\]
and the sets
\[
A_k = \{(y,t);\ y\in \U_k,\
t\in (0, \len(Y_k(y)))\},
\quad
\Omega_k = \Psi_k(A_k)\,.
\]
By construction, it is readily seen that
$q_k = 0$ outside $\Omega_k$ for every $k$,
and that $\sum_{k=1}^N q_k(x) = 1$
for every $x\in\Omega\setminus\overline{\Sigma}$,
hence, from Corollary~\ref{c:mn},
for almost every $x\in\Omega$.
The first equality in (\ref{f:chvar})
is obtained
using, for every $k=1,\ldots,N$, the change of variables
$\Psi_k$ on $A_k$,
Fubini's theorem,
and the area formula
(see \cite[\S3.3.3 and \S3.3.4]{EG})
as follows:
\[
\begin{split}
& \int_{\Omega} h(x)\, dx  =
\sum_{k=1}^N \int_{\Omega} q_k(x) h(x)\, dx
= \sum_{k=1}^N \int_{\Omega_k} q_k(x) h(x)\, dx
\\ & =
\sum_{k=1}^N
\int_{A_k} q_k(\Psi_k(y,t))\, h(\Psi_k(y,t))\, \det D\Psi_k(y,t)\, dt\, dy
\\ & =
\sum_{k=1}^N
\int_{\U_k}\left[
\int_0^{\len(Y_k(y))} q_k(Y_k(y))\, h(\Psi_k(y,t))\,
\frac{\det D\Psi_k(y,t)}{\sqrt{g^k(y)}}\, dt
\right] \sqrt{g^k(y)}\, dy
\\ & =
\sum_{k=1}^N
\int_{Y_k(\U_k)} q_k(x) \left[
\int_0^{\len(x)} h(\Phi(x,t))\,
\gauge(\nor(x))\,D(x, t)\, dt
\right]\, d\haus(x)
\\ & =
\int_{\partial\Omega} \left[
\int_0^{\len(x)} h(\Phi(x,t))\,
\gauge(\nor(x))\,D(x, t)\, dt
\right]\, d\haus(x)
\end{split}
\]
where $D(x, t) := \det(I_{n-1}-t\, \bw(x))$.
The second equality in (\ref{f:chvar}) follows upon observing that
$\curvg_1(x),\ldots,\curvg_{n-1}(x)$
are the eigenvalues of $\bw(x)$ for every $x\in\partial\Omega$.
\end{proof}

Let us consider the following system of PDEs of Monge-Kantorovich type:
\begin{equation}\label{f:syst1}
\begin{cases}
-\dive(v\, D\gauge(Du)) = f
&\textrm{in $\Omega$},\\
\gauge(Du)\leq 1
&\textrm{in $\Omega$},\\
\gauge(Du) = 1
&\textrm{in $\{v>0\}$},
\end{cases}
\end{equation}
where the source $f\geq 0$ is a continuous function in $\Omega$,
complemented with the conditions
\begin{equation}\label{f:syst2}
\begin{cases}
u\geq 0,\
v\geq 0
&\textrm{in $\Omega$},\\
u=0
&\textrm{on $\partial\Omega$}.
\end{cases}
\end{equation}
The first equation in (\ref{f:syst1}) has to be
understood in the sense of distributions, whereas
$u$ is a viscosity solution to
$\gauge(Du) = 1$ in the set $\{v>0\}$.
We look for a solution $(u,v)$ to
(\ref{f:syst1})--(\ref{f:syst2}) in the class of continuous
and non-negative functions.

Since $\dist$ is a viscosity solution of $\gauge(D u) = 1$
in $\Omega$, $u = 0$ on $\partial\Omega$,
it is clear that $\dist$ is a viscosity solution,
vanishing on $\partial\Omega$, of the equation
$\gauge(D\dist) = 1$ in $\{v>0\}$,
for every fixed continuous function
$v\in C(\Omega)$  .
Then it is enough to prove that the equation
\begin{equation}\label{f:PDE}
-\dive(v(x)\, D\gauge(D\dist(x))) = f(x)
\qquad\textrm{in $\Omega$},\\
\end{equation}
has a continuous non-negative solution.
More precisely,
we are interested in finding
a function $v\in C(\Omega)$
satisfying
\begin{equation}\label{f:weak}
\int_{\Omega}
v(x)\, \pscal{D\gauge(D\dist(x))}{D\varphi(x)}\, dx
= \int_{\Omega} f(x)\varphi(x)\, dx
\end{equation}
for every $\varphi$
belonging to the set
$C^{\infty}_c(\Omega)$
of functions of class $C^{\infty}(\Omega)$
with compact support in $\Omega$.

In order to write explicitly a solution $v$
of (\ref{f:PDE}),
it will be convenient
to extend the functions $\Phi$, defined on $\partial\Omega$
 in Theorem~\ref{t:chvar},
and $\curvg_i$, $i=1,\ldots,n-1$, defined on $\partial\Omega$
in Definition~\ref{d:curv},
to $\overline{\Omega}\setminus\overline{\Sigma}$
in the following way.
If $x\in\overline{\Omega}\setminus\overline{\Sigma}$
and $\proj(x) = \{x_0\}$, we set
\[
\Phi(x, t) := x_0 + (\dist(x)+t)\, D\gauge(\nor(x_0)),
\qquad
\curvg_i(x) := \curvg_i(x_0),\quad
i=1,\ldots,n-1\,.
\]
Since, from Lemma~\ref{l:pusc},
the map which associates to every point
$x\in\overline{\Omega}\setminus\overline{\Sigma}$
its unique projection on $\partial\Omega$
is continuous in $\overline{\Omega}\setminus\overline{\Sigma}$,
we conclude that
the maps $\curvg_i$, $i=1,\ldots,n-1$, which are continuous on $\partial\Omega$
(see Remark~\ref{r:eigv}),
are also continuous in $\overline{\Omega}\setminus\overline{\Sigma}$,
and the map $\Phi$ is continuous in
$(\overline{\Omega}\setminus\overline{\Sigma})\times\R$.

The remaining part of this section will be devoted to
the proof of the following theorem.

\begin{theorem}\label{t:existence}
Let $f\in C(\Omega)$, $f\geq 0$.
Then the function
\begin{equation}\label{f:vf}
\vf(x) =
\begin{cases}
\displaystyle{
\int_0^{\tau(x)} f(\Phi(x,t))
\prod_{i=1}^{n-1}\frac{1-(\dist(x)+t)\, \curvg_i(x)}%
{1-\dist(x)\, \curvg_i(x)}\, dt}
&\textrm{if $x\in{\overline{\Omega}}\setminus\overline{\Sigma}$},\\
0,
&\textrm{if $x\in\overline{\Sigma}$}\,,
\end{cases}
\end{equation}
is continuous in $\Omega$ and satisfies (\ref{f:weak})
for every $\varphi\in C^{\infty}_c(\Omega)$.
Here $\tau$ is the distance to cut locus
defined in (\ref{f:tau}).
\end{theorem}

Before proving Theorem~\ref{t:existence},
we establish some basic bounds on $\vf$.

\begin{lemma}
Let $f\in C(\Omega)$, $f\geq 0$.
Then, in any set $\Omega_{\epsilon} := \{x\in\Omega;\
\dist(x)>\epsilon\}$, $\epsilon > 0$, $\vf$ satisfies the bounds
\begin{equation}\label{f:estiv}
0\leq \vf(x)\leq
\| f\|_{C(\Omega_{\epsilon})}\,\tau(x)
\prod_{i=1}^{n-1} (1+T\, \widetilde{K}_{-})
\qquad
\forall x\in \Omega_{\epsilon}\,,
\end{equation}
where $T$ and $\widetilde{K}_{-}$
are the constants defined by
\begin{gather}
T=\max\{\tau(x);\ x\in\partial\Omega\},\label{f:T}\\
\widetilde{K}_{-} = \max\{{[\curvg_i(x)]}_{-};\
x\in\partial\Omega,\ i=1,\ldots,n-1\},
\label{f:K}
\end{gather}
being
$[a]_{-} = \max\{0, -a\}$ the negative part
of a real number $a$.
\end{lemma}

\begin{proof}
By Lemma~\ref{l:W2} and Corollary~\ref{c:closed}
we have that
$1-(\dist(x)+t)\curvg_i(x)>0$ for every
$x\in \overline{\Omega}\setminus\overline{\Sigma}$ and
$0\leq t\leq \tau(x)$. Then the function $v_f$ is well defined
and
$\vf\geq 0$ in $\Omega$.
In order to prove (\ref{f:estiv})
it is enough to observe that
\[
\frac{1-(\dist(x)+t)\, \curvg_i(x)}%
{1-\dist(x)\, \curvg_i(x)} \leq
1+t\, [\curvg_i(x)]_{-}
\]
and that $x+t\, D\gauge(D\dist(x))\in\Omega_{\epsilon}$
for every $x\in\Omega_{\epsilon}$ and $t\in [0,\tau(x)]$.
\end{proof}

\begin{proof}[Proof of Theorem~\ref{t:existence}]
The continuity of $\vf$ in ${\Omega}\setminus\overline{\Sigma}$
follows from
the continuity of the functions
$\curvg_i$, $i=1,\ldots,n$,
in $\overline{\Omega}\setminus\overline{\Sigma}$,
of the functions $f$ and $\tau$ in $\Omega$,
and of the function $\Phi$ in
$(\Omega\setminus\overline{\Sigma})\times\R$.
On the other hand, the continuity of $\vf$ on $\overline{\Sigma}$
is a consequence of (\ref{f:estiv}).

Let $\varphi\in C^{\infty}_c(\Omega)$.
From the change of variables formula (\ref{f:chvar})
we have that
\[
\begin{split}
& \int_{\Omega} f(x)\varphi(x)\, dx
\\ & =
\int_{\partial\Omega} \gauge(\nor(x))\, \left[
\int_{0}^{\len(x)}
f(\Phi(x,t))\,\varphi(\Phi(x,t))\,
\prod_{i=1}^{n-1} (1-t\, \curvg_i(x))
\, dt\right]\,d\haus(x)\,.
\end{split}
\]
Let us compute the term in brackets, integrating by parts,
and taking into account
that
$\varphi(\Phi(x,0)) = \varphi(x) = 0$
for every $x\in\partial\Omega$.
Setting $\Phi=\Phi(x,t)$, we get
\begin{equation}\label{f:finfin}
\begin{split}
&\int_{0}^{\len(x)}
f(\Phi)\,\varphi(\Phi)\,
\prod_{i=1}^{n-1} (1-t\, \curvg_i(x))
\, dt  \\
& =
\int_{0}^{\len(x)}
\pscal{D\varphi(\Phi)}{D\gauge(\nor(x))}
\int_{t}^{\len(x)}
f(\Phi(x,s))\,
\prod_{i=1}^{n-1} (1-s\, \curvg_i(x))\, ds
\, dt\,.
\end{split}
\end{equation}
Noticing that
\[
\vf(\Phi(x,t)) =
\int_{t}^{\len(x)}
f(\Phi(x,s))\,
\prod_{i=1}^{n-1}\frac{1-s\, \curvg_i(x)}%
{1-t\, \curvg_i(x)}\, ds
\qquad
x\in\partial\Omega,\
t\in [0,\len(x))\,,
\]
we obtain
\[
\begin{split}
& \int_{\Omega} f(x)\varphi(x)\, dx
\\ & =
\int_{\partial \Omega}\gauge(\nor)\left[
\int_{0}^{\len(x)}
\pscal{D\varphi(\Phi)}{D\gauge(\nor)}
\vf(\Phi)
\prod_{i=1}^{n-1} (1-t\, \curvg_i)\, dt\right]\, d\haus(x)\,.
\end{split}
\]
Finally, again from (\ref{f:chvar}),
\[
\int_{\Omega} f(x)\varphi(x)\, dx =
\int_{\Omega}
\vf(x)\, \pscal{D\gauge(D\dist(x))}{D\varphi(x)}\, dx\,,
\]
for every $\varphi\in C^{\infty}_c(\Omega)$.
\end{proof}

%
\bibliographystyle{amsplain}

\end{document}